**Finding optimal solutions for vehicle routing problem with pickup and delivery services with time windows: A dynamic programming approach based on state–space–time network representations**

Monirehalsadat Mahmoudi
mmahmoudi@asu.edu

Xuesong Zhou* (Corresponding Author)
xzhou74@asu.edu

School of Sustainable Engineering and the Built Environment, Arizona State University, Tempe, AZ 85281, USA

*Corresponding author. Tel.: +0014809655827.

**Abstract**

Optimization of on-demand transportation systems and ride-sharing services involves solving a class of complex vehicle routing problems with pickup and delivery with time windows (VRPPDTW). This paper first proposes a new time-discretized multi-commodity network flow model for the VRPPDTW based on the integration of vehicles' carrying states within space–time transportation networks, so as to allow a joint optimization of passenger-to-vehicle assignment and turn-by-turn routing in congested transportation networks. Our three-dimensional state–space–time network construct is able to comprehensively enumerate possible transportation states at any given time along vehicle space–time paths, and further allows a forward dynamic programming solution algorithm to solve the single vehicle VRPPDTW problem. By utilizing a Lagrangian relaxation approach, the primal multi-vehicle routing problem is decomposed to a sequence of single vehicle routing sub-problems, with Lagrangian multipliers for individual passengers' requests being updated by sub-gradient-based algorithms. We further discuss a number of search space reduction strategies and test our algorithms, implemented through a specialized program in C++, on medium-scale and large-scale transportation networks, namely the Chicago sketch and Phoenix regional networks.

**Keywords:** Vehicle routing problem with pickup and delivery with time windows; Lagrangian relaxation; Time-dependent least-cost path problem; Forward dynamic programming; Ride-sharing service optimization.

# 1 Introduction

As population and personal travel activities continue to increase, traffic congestion has remained as one of the major concerns for transportation system agencies with tight resource constraints. The next generation of transportation system initiatives aims to integrate various demand management strategies and traffic control measures to actively achieve mobility, environment, and sustainability goals. A number of approaches hold promises of reducing the undesirable effects of traffic congestion due to driving-alone trips, to name a few, demand-responsive transit services, dynamic ride-sharing, and intermodal traffic corridor management.

The optimized and coordinated ride-sharing services provided by transportation network companies (TNC) can efficiently utilize limited vehicle and driver resources while satisfying time-sensitive origin-to-destination transportation service requests. In a city with numerous travelers with different purposes, each traveler has his own traveling schedule. Instead of using his own car, the traveler can (by the aid of ride-sharing) bid and call a car just a few minutes before leaving his origin, or preschedule a car a day prior to his departure. The on-demand transportation system provides a traveler with a short waiting time even if he resides in a high-demand area. Currently, several real-time ride-sharing or, more precisely, app-based transportation network and taxi companies, such as Uber and Lyft are serving passengers in many metropolitan areas. In the long run, a fully automated and optimized ride-sharing approach is expected to handle very complex transportation supply-to-demand assignment tasks and offer a long list of benefits for transportation road users and TNC operators. These benefits might include reducing driver stress and driving cost, improving mobility for non-drivers, increasing safety and fuel efficiency, and decreasing road congestion as well as reducing overall societal energy use and pollution.

The ride-sharing problem can be mathematically modeled by one of the well-known optimization problems which is the vehicle routing problem with pickup and delivery (VRPPD). In this paper, in order to improve the solution quality and computational efficiency of on-demand transportation systems and dynamic ride-sharing services, especially for large-scale real-world transportation networks, we propose a new mathematical programming model for the vehicle routing problem with pickup and delivery with time windows (VRPPDTW) that can fully recognize time-dependent link travel time caused by traffic congestion at different times of day. Based on the Lagrangian relaxation solution framework, we further present a holistic optimization approach for matching passengers' requests to transportation service providers, synchronizing transportation vehicle routing, and determining request pricing (e.g. through Lagrangian multipliers) for balancing transportation demand satisfaction and resource needs on urban networks.

# 2 Literature review and research motivations

The vehicle routing problem with pickup and delivery with time windows (VRPPDTW) or simply, pickup and delivery problem with time windows (PDPTW), is a generalized version of the vehicle routing problem with time windows (VRPTW), in which each transportation request is a combination of pickup at the origin node and drop-off at the destination node (Desaulniers et al., 2002). The PDPTW under consideration in this paper contains all constraints in the VRPTW plus added constraints in which either pickup or delivery has given time windows, and each request must be served by a single vehicle. The PDPTW may be observed as the dial-a-ride problem in the literature as well. Since the VRPTW is an NP-hard problem, the PDPTW is also NP-hard (Baldacci et al., 2011).

Several applications of the VRPPDTW have been reported in road, maritime, and air transportation environments, to name a few, Fisher et al. (1982), Bell et al. (1983), Savelsbergh and Sol (1998), Wang and Regan (2002), and Zachariadis et al. (2015) in road cargo routing and scheduling; Psaraftis et al. (1985), Fisher and Rosenwein (1989), and Christiansen (1999) in sea cargo routing and scheduling; and Solanki and Southworth (1991), Solomon et al. (1992), Rappoport et al. (1992), and Rappoport et al. (1994) in air cargo routing and scheduling. Further applications of the VRPPDTW can be found in transportation of elderly or handicapped people (Jaw et al., 1986; Alfa, 1986; Ioachim et al., 1995; and Toth and Vigo, 1997), school bus routing and scheduling (Swersey and Ballard, 1983; and Bramel and Simchi-Levi, 1995), and ride-sharing (Hosni et al., 2014; and Wang et al., 2015). Recently, Furuhata et al. (2013) offers an excellent review and provides a systematic classification of emerging ride-sharing systems.

Although clustering algorithms (Cullen et al., 1981; Bodin and Sexton, 1986; Dumas et al., 1989; Desrosiers et al., 1991; and Ioachim et al., 1995), meta-heuristics (Gendreau et al., 1998; Toth and Vigo, 1997; and Paquette et al., 2013), neural networks (Shen et al., 1995), and some heuristics such as double-horizon based heuristics (Mitrovic-Minic et al., 2004) and regret insertion heuristics (Diana and Dessouky, 2004) have been shown to be efficient in solving a particular size of PDPTW, in general, finding the exact solution via optimization approaches has still remained theoretically and computationally challenging. Focusing on the PDPTW for a single vehicle, Psaraftis (1980) presented an exact backward dynamic programming (DP) solution algorithm to minimize a weighted combination of the total service time and the total waiting time for all customers with $O(n^2 3^n)$ complexity. Psaraftis (1983) further modified the algorithm to a forward DP approach. Sexton and Bodin (1985a, b) decomposed the single vehicle PDPTW to a routing problem and a scheduling sub-problem, and then they applied Benders' decomposition for both master problem and sub-problem, independently. Based on a static network flow formulation, Desrosiers et al. (1986) proposed a forward DP algorithm for the single-vehicle PDPTW with the objective function of minimizing the total traveled distance to serve all customers. After presenting our proposed model in the later section, we will conduct a more systematical comparison between our proposed state–space–time DP framework and the classical work by Psaraftis (1983) and Desrosiers et al. (1986).

There are a number of studies focusing on the multi-vehicle pickup and delivery problem with time windows. Dumas et al. (1991) proposed an exact algorithm to the multiple vehicle PDPTW with multiple depots, where the objective is to minimize the total travel cost with capacity, time window, precedence and coupling constraints. They applied a column generation scheme with a shortest path sub-problem to solve the PDPTW, with tight vehicle capacity constraints, and a small size of requests per route. Ruland (1995) and Ruland and Rodin (1997) proposed a polyhedral approach for the vehicle routing problem with pickup and delivery. Savelsbergh and Sol (1998) proposed an algorithm for the multiple vehicle PDPTW with multiple depots to minimize the number of vehicles needed to serve all transportation requests as the primary objective function, and minimizing the total distance traveled as the secondary objective function. Their algorithm moves toward the optimal solution after solving the pricing sub-problem using heuristics. They applied their algorithm for a set of randomly generated instances. In a two-index formulation proposed by Lu and Dessouky (2004), a branch-and-cut algorithm was able to solve problem instances. Cordeau (2006) proposed a branch-and-cut algorithm based on a three-index formulation. Ropke et al. (2007) presented a branch-and-cut algorithm to minimize the total routing cost, based on a two-index formulation. Ropke and Cordeau (2009) presented a new branch-and-cut-and-price algorithm in which the lower bounds are computed by the column generation algorithm and improved by introducing different valid inequalities to the problem. Based on a set-partitioning formulation improved by additional cuts, Baldacci et al. (2011) proposed a new exact algorithm for the PDPTW with two different objective functions: the primary is minimizing the route costs, whereas the secondary is to minimize the total vehicle fixed costs first, and then minimize the total route costs.

Previous research has made a number of important contributions to this challenging problem along different formulation or solution approaches. On the other hand, there are a number of modeling and algorithmic challenges for a large-scale deployment of a vehicle routing and scheduling algorithm, especially for regional networks with various road capacity and traffic delay constraints on freeway bottlenecks and signal timing on urban streets. A few previous research directly considers the underlying transportation network with time of day traffic congestion (Kok et al. 2012, Gromicho et al. 2012) and has defined the PDPTW on a directed graph containing customers' origin and destination locations connected by some links which are representative of the shortest distance or least travel time routes between origin–destination pairs. That is, with each link, there are associated routing cost and travel time between the two service nodes. Unlike the existing offline network for the PDPTW in which each link has a fixed routing cost (travel time), our research particularly examines the PDPTW on real-world transportation networks containing a transportation node-link structure in which routing cost (travel time) along each link may vary over the time.

In order to consider many relevant practical aspects, such as waiting costs at different locations, we utilize space–time scheme (Hägerstrand, 1970; Miller, 1991, Ziliaskopoulos and Mahmassani, 1993) to formulate the PDPTW on state–space–time transportation networks. The constructed networks are able to conveniently represent the complex pickup and delivery time windows without adding the extra constraints typically needed for the classical PDPTW formulation (e.g. Cordeau, 2006). The introduced state–space–time networks also enable us to embed computationally efficient dynamic programming algorithms for solving the PDPTW without relying on off-the-shelf optimization solvers. Even though the solution space created by our formulation has multiple dimensions and accordingly large in its sizes, the readily available large amount of computer memory in modern workstations can easily accommodate the multi-dimensional solution vectors utilized in our application. Our fully customized solution algorithms, implemented in an advanced programming language such as C++, hold the promise of tackling large-sized regional transportation network instances. To address the multi-vehicle assignment requirement, we relax the transportation request satisfaction constraints into the objective function and utilize the related Lagrangian relaxation (LR) solution framework to decompose the primal problem to a sequence of time-dependent least-cost-path sub-problems.

In our proposed solution approach, we aim to incorporate several lines of pioneering efforts in different directions. Specifically, we (1) reformulate the VRPPDTW as a time-discretized, multi-dimensional, multi-commodity flow model with linear objective function and constraints, (2) extend the static DP formulation to a fully time-dependent DP framework for single-vehicle VRPPDTW problems, and (3) develop a LR solution procedure to decompose the multi-vehicle scheduling problem to a sequence of single-vehicle problems and further nicely integrate the demand satisfaction multipliers within the proposed state–space–time network.

Based on the Lagrangian relaxation solution framework, we further present a holistic optimization approach for matching passengers' requests to transportation service providers, synchronizing transportation vehicle routing, and determining request pricing (e.g. through Lagrangian multipliers) for balancing transportation demand satisfaction and resource needs on urban networks.

The rest of the paper is organized as follows. Section 3 contains a precise mathematical description of the PDPTW in the state–space–time networks. In section 4, we present our new integer programming model for the PDPTW followed by a comprehensive comparison between Cordeau's model and our model. Then, we will show how the main problem is decomposed to an easy-to-solve problem by the Lagrangian relaxation algorithm in section 5. Section 6 provides computational results of the six-node transportation network, followed by the Chicago sketch and Phoenix regional networks to demonstrate the computational efficiency and solution optimality of our developed algorithm coded by C++. After large-scale network experiments, we conclude the paper in section 7 with discussions on possible extensions.

## 3 Problem statement based on state–space–time network representation

In this section, we first introduce our new mathematical model for the PDPTW. This is followed by a comprehensive comparison between our proposed model and the three-index formulation of Cordeau (2006) for the PDPTW, presented in Appendix A, for the demand node-oriented network.

### 3.1 Description of the PDPTW in state–space–time networks

We formulate the PDPTW on a transportation network, represented by a directed graph and denoted as $G(N, A)$, where $N$ is the set of nodes (e.g. intersections or freeway merge points) and $A$ is the set of links with different link types such as freeway segments, arterial streets, and ramps. As shown in Table 1, each directed link $(i, j)$ has time-dependent travel time $TT(i, j, t)$ from node $i$ to node $j$ starting at time $t$. Every passenger $p$ has a preferred time window for departure from his origin, $[a_p, b_p]$, and a desired time window for arrival at his destination, $[a'_p, b'_p]$, where $a_p$, $b_p$, $a'_p$, and $b'_p$ are passenger $p$'s earliest preferred departure time from his origin, latest preferred departure time from his origin, earliest preferred arrival time at his destination, and latest preferred arrival time at his destination, respectively. Each vehicle $v$ also has the earliest departure time from its starting depot, $e_v$, and the latest arrival time at its ending depot, $l_v$. In the PDPTW, passengers may share their trip with each other; in other words, every vehicle $v$, considering its capacity $Cap_v$ and the total routing cost, may serve as many passengers as possible provided that passenger $p$ is picked up and dropped-off in his preferred time windows, $[a_p, b_p]$ and $[a'_p, b'_p]$, respectively.

Each transportation node has the potential to be the spot for picking up or dropping off a passenger. Likewise, a vehicle's depot might be located at any node in the transportation network. To distinguish regular transportation nodes from passengers' and vehicles' origin and destination, we add a single dummy node $o'_v$ for vehicle $v$'s origin depot and a single dummy node $d'_v$ for vehicle $v$'s destination depot. Similarly, we can also add dummy nodes $o_p$ and $d_p$ for passenger $p$. Each added dummy node is only connected to its corresponding physical transportation node by a link. The travel time on this link can be interpreted as the service time if the added dummy node is related to a passenger's origin or destination, and as preparation time if it is related to a vehicle's starting or ending depot. Table 1 lists the notations for the key sets, indices and parameters in the PDPTW.

**Table 1**
Sets, indices and parameters in the PDPTW.

| Symbol | Definition |
|---|---|
| $V$ | Set of physical vehicles |
| $V^*$ | Set of virtual vehicles |
| $P$ | Set of passengers |
| $N$ | Set of physical transportation nodes in the physical traffic network based on geographical location |
| $W$ | Set of possible passenger carrying states |
| $v$ | Vehicle index |
| $v_p^*$ | Index of virtual vehicle exclusively dedicated for passenger $p$ |
| $p$ | Passenger index |
| $w$ | Passenger carrying state index |
| $(i,j)$ | Index of physical link between adjacent nodes $i$ and $j$ |
| $TT(i,j,t)$ | Link travel time from node $i$ to node $j$ starting at time $t$ |
| $Cap_v$ | Maximum capacity of vehicle $v$ |
| $a_p$ | Earliest departure time from passenger $p$'s origin |
| $b_p$ | Latest departure time from passenger $p$'s origin |
| $a'_p$ | Earliest arrival time at passenger $p$'s destination |
| $b'_p$ | Latest arrival time at passenger $p$'s destination |
| $[a_p, b_p]$ | Departure time window for passenger $p$'s origin |
| $[a'_p, b'_p]$ | Arrival time window for passenger $p$'s destination |
| $o'_v$ | Dummy node for vehicle $v$'s origin |
| $d'_v$ | Dummy node for vehicle $v$'s destination |
| $e_v$ | Vehicle $v$'s earliest departure time from the origin depot |
| $l_v$ | Vehicle $v$'s latest arrival time to the destination depot |
| $o_p$ | Dummy node for passenger $p$'s origin (pickup node for passenger $p$) |
| $d_p$ | Dummy node for passenger $p$'s destination (delivery node for passenger $p$) |

We now use an illustrative example to demonstrate key modeling features of constructed networks. Consider a physical transportation network consisting of six nodes presented in Fig. 1. Each link in this network is associated with time-dependent travel time $TT(i,j,t)$. Without loss of generality, the number written on each link denotes the time-invariant travel time $TT(i,j)$ in terms of minutes. Suppose two requests with two origin–destination pairs should be served. For simplicity, it is assumed that both passengers have the same origin (node 2) and the same drop-off node (node 3). There is only one vehicle available for serving. Moreover, it is assumed that the vehicle starts its route from node 4 and ends it at node 1. Passenger 1 should be picked up from dummy node $o_1$ in time window $[4,7]$ and dropped off at dummy node $d_1$ in time window $[11,14]$, while Passenger 2 should be picked up from dummy node $o_2$ in time window $[8,10]$ and dropped off at dummy node $d_2$ in time window $[13,16]$. Vehicle 1 also has the earliest departure time from its starting depot, $t = 1$, and the latest arrival time at its ending depot, $t = 20$.

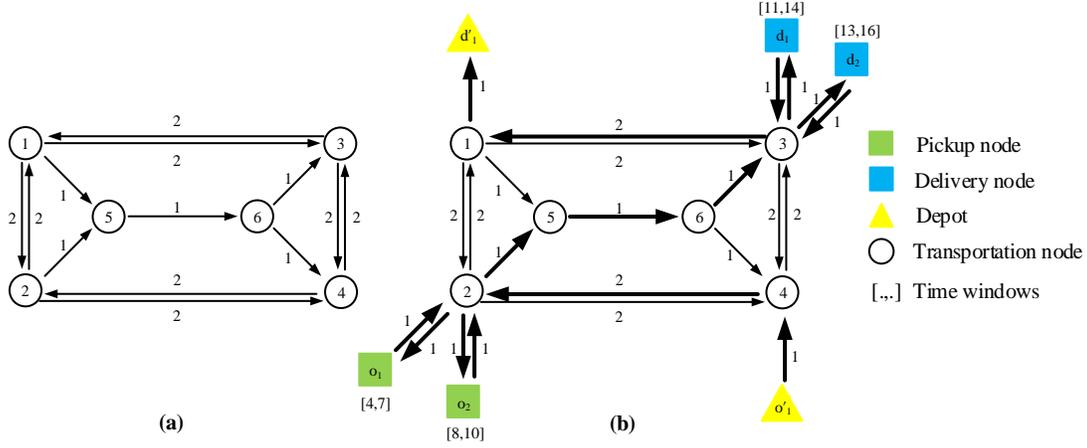

**Fig. 1.** (a) Six-node transportation network; (b) transportation network with the corresponding dummy nodes.

Note that the shortest path with node sequence $(o'_1, 4, 2, o_1, 2, o_2, 2, 5, 6, 3, d_1, 3, d_2, 3, 1, d'_1)$ from vehicle 1's origin to its ending depot is shown by bold arrows when it serves both passenger 1 and 2. To construct a state–space–time network, the time horizon is discretized into a series of time intervals with the same time length. Without loss of generality, we assume that a unit of time has 1 min length. Interested readers are referred to Yang and Zhou (2014) on details about how to construct a space–time network. To avoid more complexity in the vehicle's space–time network illustrated in Fig. 2, only those arcs constituting the shortest paths from vehicle 1's origin to its destination are demonstrated. Our formulation has a set of precise rules to allow or restrict the vehicle waiting behavior in the constructed space–time network, depending on the type of nodes and the associated time window. First, vehicle $v$ may wait at its own origin or destination depot or at any other physical transportation nodes. If a vehicle arrives at passenger $p$'s origin node before time $a_p$, it must wait at the related physical node until the service is allowed to begin. Moreover, we assume that a vehicle is not allowed to stop at passenger $p$'s dummy origin node after time $b_p$. Similarly, if a vehicle arrives at passenger $p$'s destination node before time $a'_p$, it must wait until it is allowed to drop-off passenger $p$, and vehicle $v$ is not allowed to stop at passenger $p$'s dummy destination node after time $b'_p$.

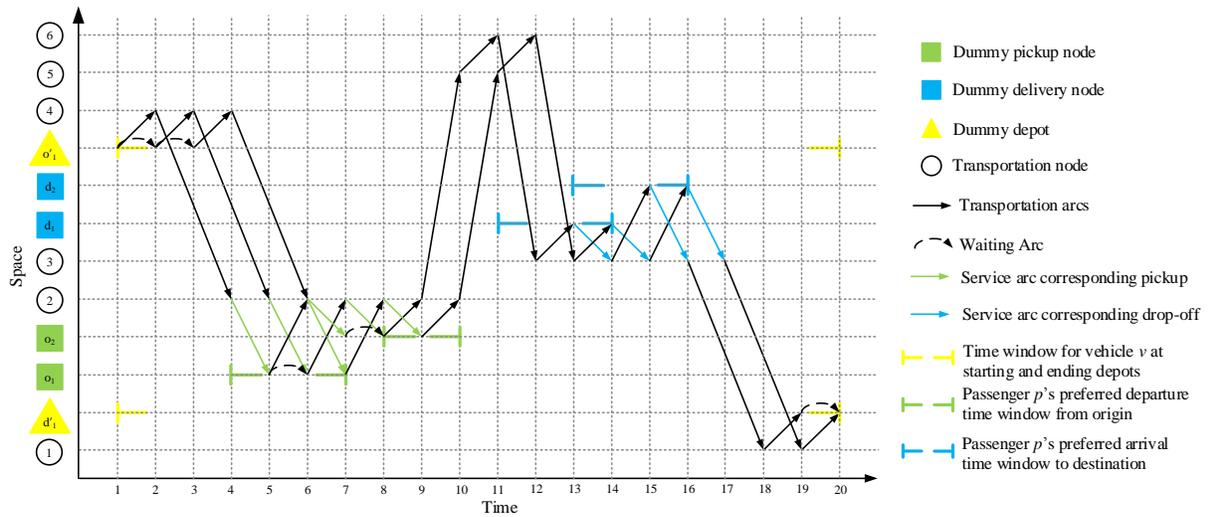

**Fig. 2.** Shortest paths with node sequence $(o'_1, 4, 2, o_1, 2, o_2, 2, 5, 6, 3, d_1, 3, d_2, 3, 1, d'_1)$ in vehicle 1's space–time network.

In the problem under consideration, we assume all passengers' desired departure and arrival time windows are feasible. However, it is quite possible that some passenger transporting requests could not be satisfied at all since the total number of physically available vehicles in the ride-sharing company or organization is not enough to satisfy all the demands. To avoid infeasibility for the constructed optimization problem, we define a virtual vehicle for each passenger exclusively. We assume that both starting and ending depots of virtual vehicle $v_p^*$ are located exactly where passenger $p$ is going to be picked up. By doing so, there is no cost incurred if the virtual vehicle is not needed to carry the related passenger, and in this case the virtual vehicle simply waits at its own depot. On the other hand, if the virtual vehicle is needed to perform the service to ensure there is a feasible solution, then virtual vehicle $v_p^*$ starts its route from its starting depot, picks up passenger $p$, delivers him to his destination, and then comes back to its ending depot. Fig. 3 shows the shortest paths with node sequence $(o'_{1^*}, 2, o_1, 2, 5, 6, 3, d_1, 3, 1, 2, d'_{1^*})$ in vehicle $v_1^*$'s space–time network.

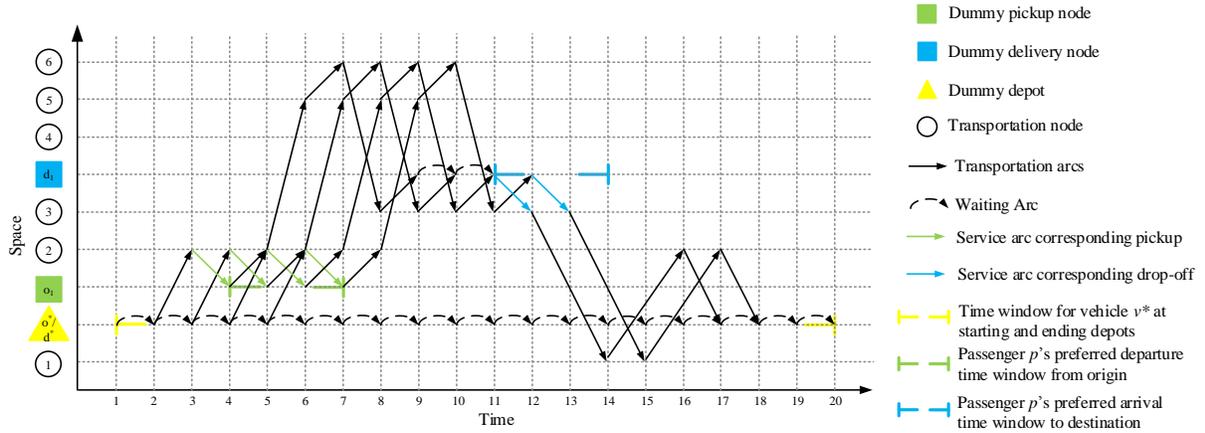

**Fig. 3.** Shortest paths with node sequence $(o'_{1^*}, 2, o_1, 2, 5, 6, 3, d_1, 3, 1, 2, d'_{1^*})$ in vehicle $v_1^*$'s space–time network.

### 3.2 Representing the state of system and calculating the number of states

In the context of dynamic programming, we need to decompose the complex VRP structure into a sequence of overlapping stage-by-stage sub-problems in a recursive manner. For each stage of the optimization problem, we need to define the state of the process so that the state of the system with $n$ stages to go can fully summarize all relevant information of the system for future decision-making purposes no matter how the process has reached the current stage $n$. In our pickup and delivery problem, in each vehicle's network, the given time index $t$ acts as the stage, and the state of the system is jointly defined by two indexes: node index $i$ and the passenger carrying state index $w$. The latter passenger carrying state $w$ can be also represented as a vector with $|P|$ number of elements $[\pi_1, \pi_2, ..., \pi_p, ..., \pi_P]$, where $\pi_p$ equals 1 or 0 and denotes the status of passenger $p$ whether he is riding the vehicle or not. To facilitate the descriptions of the state transition, we introduce the following equivalent notation system for passenger carrying states: if a vehicle carries passenger $p$, the $p^{\text{th}}$ element of the state $w$ is filled with passenger $p$'s id; otherwise, it is filled with a dash sign, as illustrated in Table 2.

**Table 2**
Binary representation and equivalent character-based representation for passenger carrying states.

| Binary representation | Equivalent character-based representation |
|---|---|
| [0,0,0] | [_ _ _] |
| [1,0,0] | [$p_1$ _ _] |
| [0,1,1] | [_ $p_2$ $p_3$] |

Without loss of generality, for a typical off-line vehicle routing problem, the initial and ending states of the vehicles are assumed to be empty, corresponding to the state $[\_\_\_]$. For an on-line dynamic vehicle dispatching application, one can define the starting passenger carrying state to indicate the existing passengers riding the vehicle, for example, $[p_1 \_\_]$ if passenger 1 is being served currently. We use an illustrative example to demonstrate the concept of a passenger's carrying state clearly. Suppose three requests with three different origin–destination pairs should be served. There is only one vehicle available for serving and let's assume that the vehicle can carry up to two passengers at the same time. We can enumerate all different carrying states for the vehicle. The first state is the state in which the vehicle does not carry any passenger $[\_\_\_]$. There are $C_1^3$ number of possible carrying states in which the vehicle only carries one passenger at time $t$: $[p_1 \_\_]$, $[\_ p_2 \_]$, and $[\_\_ p_3]$. Similarly, there are $C_2^3$ number of possible carrying states in which the vehicle carries two passengers at time $t$ which are $[p_1 \, p_2 \_]$, $[p_1 \_ p_3]$, and $[\_ p_2 \, p_3]$. Since the vehicle can carry up to two passengers at the same time, the state of $[p_1 \, p_2 \, p_3]$ is infeasible. Fig. 4(a) and Fig. 4(b) show shared ride state $[p_1 \, p_2 \_]$ and single-passenger-serving state $[\_ p_2 \_]$.

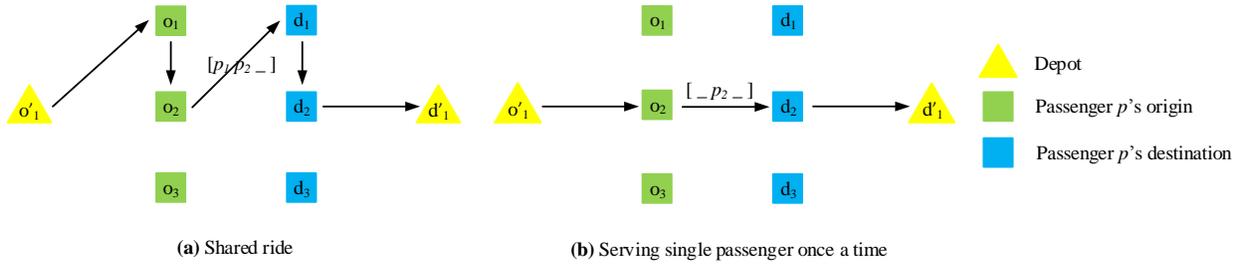

**(a)** Shared ride  **(b)** Serving single passenger once a time

**Fig. 4.** State transition path (a) Passenger carrying state $[p_1 \, p_2 \_]$; (b) Passenger carrying state $[\_ p_2 \_]$.

We are further interested in the number of feasible states, which critically determines the computational efforts of the DP-based solution algorithm. First, there is a unique state in which vehicle $v$ does not carry any passenger, which is a combinatory of $C_0^P$ for selecting 0 passengers from the collection of $P$ passengers. Similarly, there are $C_1^P$ number of possible carrying states in which vehicle $v$ only carries one passenger at a time. Likewise, there are $C_k^P$ number of possible carrying states in which vehicle $v$ carries $k$ passengers at a time. Note that $k \leq Cap_v$. Therefore, the total number of possible passenger carrying states is equal to $\sum_{k=0}^{Cap_v} C_k^P$. It should be remarked that, according to the earliest departure time from the origin and the latest arrival time to the destination of different passengers, some of the possible carrying states, say $[\_ p_2 \, p_3]$, might be infeasible as there is insufficient transportation time to pick up those two passengers together while satisfying their time window constraints.

Consider the following example, where passenger 1 should be picked up in time window [4,7] and delivered in time window [9,12], whereas passenger 3's preferred time windows for being picked up and delivered are [20,24] and [25,29], respectively. So, it is obvious that passenger 1 and 3 cannot share their ride with each other and be transported at the same time by the same vehicle. Therefore, state $[p_1 \_ p_3]$ is definitely infeasible in this example. We will further explain how to reduce the search region by defining some rational rules and simple heuristics in section 5.3.

**3.3 State transition associated with pickup and delivery links**

Each vehicle starts its trip from the empty state in which the vehicle does not carry any passengers. We call this state as the initial state ($w_0$). Each vertex in the constructed state–space–time network is recognized by a triplet of three different indexes: node index $i$, time interval index $t$, and passenger carrying state index $w$. In the space–time transportation network construct, we can identify a traveling arc $(i,j,t,s)$ starting from node $i$ at time $t$ arriving at node $j$ at time $s$. Accordingly, in the state–space–time network, each vertex $(i,t,w)$ is connected to vertex $(j,s,w')$ through arc $(i,j,t,s,w,w')$. To find all feasible combinations of passenger carrying state transition $(w,w')$ on an arc, it is sufficient to follow these rules:

Rule 1. On a pickup link (with the passenger origin dummy node as the downstream node), vehicle $v$ picks up passenger $p$, so $\pi_p$ is changed from 0 to 1, or equivalently, the $p^{th}$ element of the corresponding states should be changed from a dash sign to passenger $p$ id.

Rule 2. On a drop-off link (with the passenger destination dummy node as the upstream node), vehicle $v$ drops off passenger $p$, so $\pi_p$ is changed from 1 to 0, and the $p^{th}$ element of the corresponding states should be changed from passenger $p$ id to a dash sign.

Rule 3. On a transportation link or links connected to vehicle dummy nodes, vehicle $v$ neither picks up nor drops off any passenger, and all elements of $w$ and $w'$ should be the same.

To find all feasible passengers state transition $(w, w')$, we need to examine all possible combinations of $w$ and $w'$. Consider a three-passenger case, in which Table 3 identifies all possible combinations of these state transitions. Note that the vehicle can carry up to two passengers at the same time in this example. The empty cells indicate impossible state transitions in the constructed space–time network with dedicated dummy nodes. The corresponding possible passenger carrying state transitions (pickup or drop-off) are illustrated in one graph in Fig. 5. Fig. 6 represents the projection on state–space network for the example presented in section 3.1.

**Table 3**
All possible combinations of passenger carrying states.

| $w$ \ $w'$ | [_ _ _] | [$p_1$ _ _] | [_ $p_2$ _] | [_ _ $p_3$] | [$p_1$ $p_2$ _] | [$p_1$ _ $p_3$] | [_ $p_2$ $p_3$] |
|---|---|---|---|---|---|---|---|
| [_ _ _] | no change | pickup | pickup | pickup | | | |
| [$p_1$ _ _] | drop-off | no change | | | pickup | pickup | |
| [_ $p_2$ _] | drop-off | | no change | | pickup | | pickup |
| [_ _ $p_3$] | drop-off | | | no change | | pickup | pickup |
| [$p_1$ $p_2$ _] | | drop-off | drop-off | | no change | | |
| [$p_1$ _ $p_3$] | | drop-off | | drop-off | | no change | |
| [_ $p_2$ $p_3$] | | | drop-off | drop-off | | | no change |

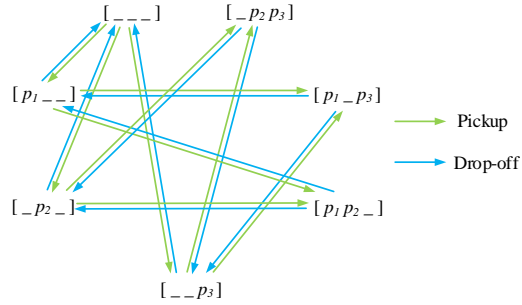

**Fig. 5.** Finite states graph showing all possible passenger carrying state transition (pickup or drop-off).

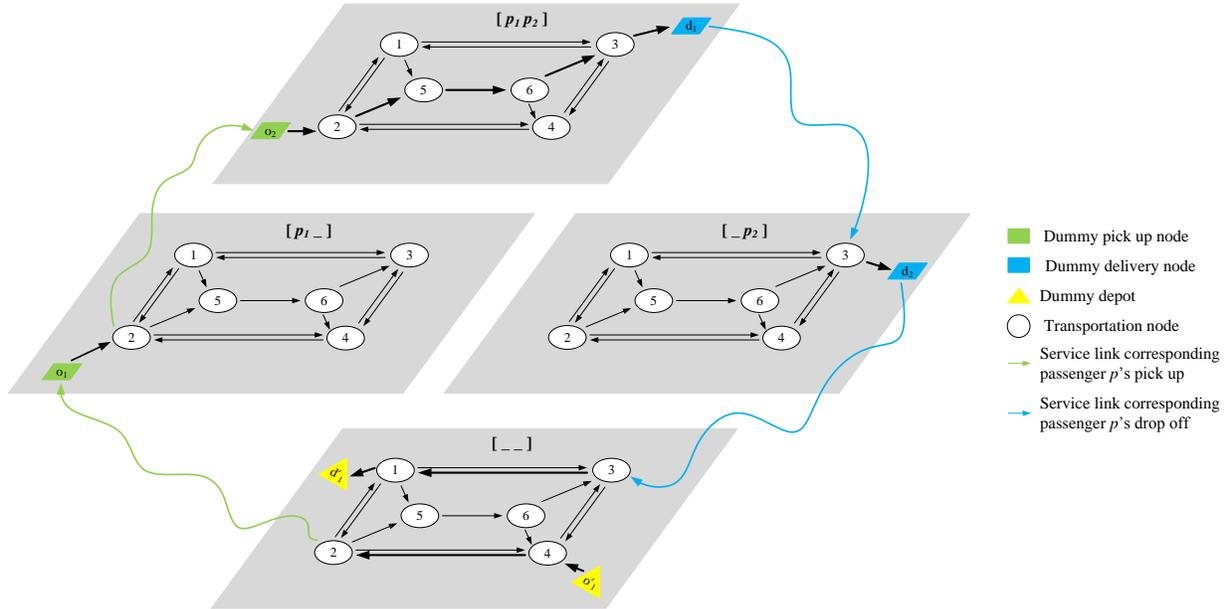

**Fig. 6.** Projection on state–space network representation for ride-sharing path (pick up passenger $p_1$ and then $p_2$).

## 4 Time-discretized multi-commodity network flow programming model

Based on the constructed state–space–time networks that can capture essential pickup and delivery time window constraints, we now start constructing a multi-commodity network flow programing model for the VRPPDTW. In this multi-dimensional network, the challenge is to systematically describe the related flow balance constraints for vehicles and request satisfaction constraints for passengers. As shown in Table 4, we use $(i, t, w)$ to represent the indices of state–space–time vertexes, and the corresponding arc index which is $(i, j, t, s, w, w')$. Let $B_v$ denote the set of state–space–time arcs in vehicle $v$'s network, which has three different types of arcs, namely, service arcs, transportation arcs and waiting arcs.

  i. All passenger carrying state transitions (i.e., pickup or drop-off) occurs only on service arcs. In other words, all incoming arcs to passengers' origin nodes (pickup arcs shown by green lines in Figures 5 and 6) and all outgoing arcs from their destination nodes (drop-off arcs shown by blue lines in Figures 5 and 6) are considered service arcs.
 ii. A link with both ends as physical nodes or vehicle dummy nodes are considered transportation arcs.
iii. Vehicles (both physical and virtual) may wait at their own origin or destination depot or at any other physical transportation nodes through waiting arcs $(i, i, t, t + 1, w, w)$ from time $t$ to time $t + 1$ with the same passenger carrying state $w$.

**Table 4**
Indexes and variables used in the time-discretized network flow model.

| Symbol | Definition |
|---|---|
| $(i, t, w), (j, s, w')$ | Indexes of state–space–time vertexes |
| $(i, j, t, s, w, w')$ | Index of a space–time-state arc indicating that one can travel from node $i$ at time $t$ with passenger carrying state $w$ to the node $j$ at time $s$ with passenger carrying state $w'$ |
| $B_v$ | Set of state–space–time arcs in vehicle $v$'s network |
| $c(v, i, j, t, s, w, w')$ | Routing cost of arc $(i, j, t, s, w, w')$ traveled by vehicle $v$ |
| $TT(v, i, j, t, s, w, w')$ | Travel time of arc $(i, j, t, s, w, w')$ traveled by vehicle $v$ |
| $\Psi_{p,v}$ | Set of pickup service arcs of passenger $p$ in vehicle $v$'s networks |
| $y(v, i, j, t, s, w, w')$ | $= 1$ if arc $(i, j, t, s, w, w')$ is used by vehicle $v$; $= 0$ otherwise |

In general, the travel time $TT(v,i,j,t,s,w,w')$ is the travel time of traversing from node $i$ at time $t$ with passenger carrying state $w$ to node $j$ at time $s$ with passenger carrying state $w'$ by vehicle $v$. As we mentioned before, travel time for service arcs can be interpreted as the service time needed to pick up or drop-off a passenger, and as the preparation time if the arc is related to a vehicle's starting or ending depot. In addition, the travel time of the waiting arcs is assumed to be a unit of time.

The routing cost $c(v,i,j,t,s,w,w')$ for an arc can be defined as follows. The routing cost of a transportation arc is defined as a ratio of its travel time. For the physical vehicle, this ratio is basically the total transportation cost per hour when the vehicle is traveling, which may include the fuel, maintenance, depreciation, insurance costs, and more importantly, the cost of hiring a full-time or part-time driver. Let's assume that, in total, the transportation by a physical vehicle costs $x$ dollars per hour. Since passengers should be served by physical vehicles by default and virtual vehicles serve passengers only if there is no available physical vehicle to satisfy their demand, we impose a quite expensive transportation cost per hour for virtual vehicles, let's say $2x$ dollars per hour. The routing cost of the service arcs are defined similarly to the routing cost of the transportation arcs. The routing cost of a waiting arc is also defined as a ratio of its waiting time. However, this ratio is basically the total cost of the physical vehicle $v$ per hour when the driver has turned off the vehicle and is waiting at a node, which may only include the cost of hiring a full-time or part-time driver. Let's assume that, in total, waiting at a node by a physical vehicle costs $y$ dollars per hour, with a typical relationship of waiting cost < transportation cost per hour, i.e., $y < x$. We assume that waiting at the origin and destination depot for a physical vehicle has no charge for the service provider in order to encourage a vehicle to reduce the total transportation time, if possible. Moreover, for virtual vehicles, the waiting cost is always equal to zero to allow a virtual vehicle be totally idle at its own depot.

The model uses binary variables $y(v,i,j,t,s,w,w')$ equal to 1 if and only if state–space–time arc $(i,j,t,s,w,w')$ is used by vehicle $v$. Without loss of generality, we assume that a vehicle does not carry any passenger when it departs from its origin depot or arrives to its destination depot, which correspond to the passenger carrying state at node $(i = o'_v, t = e_v)$ and $(j = d'_v, s = l_v)$ as an empty state denoted by $w_0$. Note that, since passenger carrying state transitions only occur through service arcs, $w = w' = w_0$ for $y(v, o'_v, j, e_v, s, w, w')$ and $y(v, i, d'_v, t, l_v, w, w')$.

Note that each vehicle must end its route at the destination depot with the empty passenger carrying state. Therefore, if vehicle $v$ picks up passenger $p$ from his origin, to maintain the flow balance constraints, the vehicle must drop-off the passenger at his destination node so that the vehicle comes back to its ending depot with the empty passenger carrying state. As a result, constraints corresponding to the passengers' drop-off request is redundant and it does not need to enter into the following formulation. After constructing the state–space–time transportation network for each vehicle, the PDPTW can be formulated as follows:

$Min\ Z = \sum_{v \in (V \cup V^*)} \sum_{(i,j,t,s,w,w') \in B_v} c(v,i,j,t,s,w,w') y(v,i,j,t,s,w,w')$ (1)

s.t.

Flow balance constraints at vehicle $v$'s origin vertex

$\sum_{(i,j,t,s,w,w') \in B_v} y(v,i,j,t,s,w,w') = 1 \qquad i = o'_v, t = e_v, w = w' = w_0, \forall v \in (V \cup V^*)$ (2)

Flow balance constraint at vehicle $v$'s destination vertex

$\sum_{(i,j,t,s,w,w') \in B_v} y(v,i,j,t,s,w,w') = 1 \qquad j = d'_v, s = l_v, w = w' = w_0, \forall v \in (V \cup V^*)$ (3)

Flow balance constraint at intermediate vertex

$\sum_{(j,s,w'')} y(v,i,j,t,s,w,w'') - \sum_{(j',s',w')} y(v,j',i,s',t,w',w) = 0 \quad (i,t,w) \notin \{(o'_v, e_v, w_0), (d'_v, l_v, w_0)\}, \forall v \in (V \cup V^*)$ (4)

Passenger $p$'s pickup request constraint

$\sum_{v \in (V \cup V^*)} \sum_{(i,j,t,s,w,w') \in \Psi_{p,v}} y(v,i,j,t,s,w,w') = 1 \qquad \forall p \in P$ (5)

Binary definitional constraint

$$y(v,i,j,t,s,w,w') \in \{0,1\} \qquad \forall (i,j,t,s,w,w') \in B_v, \forall v \in (V \cup V^*) \quad (6)$$

The objective function (1) minimizes the total routing cost. Constraints (2) to (4) ensure flow balance on every vertex in vehicle $v$'s state–space–time transportation network. Constraints (5) express that each passenger is picked up exactly once by a vehicle (either physical or virtual). Constraint (6) defines that the decision variables are binary.

The three-index formulation of Cordeau (2006) for the PDPTW in the origin–destination network is presented in Appendix A. Table 5 shows that our proposed model encompasses all constraints used in Cordeau's model.

**Table 5**
An analogy between Cordeau's model and our model for the PDPTW.

| Cordeau (2006) | Our model |
| --- | --- |
| **Three-index variables $x_{ij}^v$ for vehicle $v$ on link $(i,j)$** | Seven-index variable $y(v,i,j,t,s,w,w')$ for vehicle v on arc $(i,j,t,s,w,w')$. |
| (A.1) minimizes the total routing cost. | (1) minimizes the total routing cost. |
| (A.2) guarantees that each passenger is picked up. | (5) guarantees that each passenger is picked up. |
| (A.2) and (A.3) ensure that each passenger's origin and destination are visited exactly once by the same vehicle. | (2) to (5) ensure that the same vehicle $v$ transports passenger $p$ from his origin to his destination. |
| (A.4) expresses that each vehicle starts its route from the origin depot. | (2) expresses that each vehicle starts its route from the origin depot. |
| (A.5) ensures the flow balance on each node. | (2) to (4) ensure flow balance on every vertex in vehicle $v$'s network. |
| (A.6) expresses that each vehicle ends its route to the destination depot. | (3) expresses that each vehicle ends its route to the destination depot. |
| (A.7) ensures the validity of the time variables. | The essence of state–space–time networks ensures the time variables are calculated correctly through arc $(i,j,t,s,w,w')$, where arrival time $s = t + TT(i,j,t)$. |
| (A.8) ensures the validity of the load variables. | The structure of state–space–time networks ensures that each vehicle transports a number of passengers up to its capacity at a time, in terms of feasible states $(w,w')$. |
| (A.9) defines each passenger's ride time. | Employing state–space–time networks defines each passenger's ride time. |
| (A.10) imposes the maximal duration of each route. | Vehicle $v$'s network is constructed subject to time window $[e_v, l_v]$. |
| (A.11) imposes time windows constraints. | Passenger $p$'s network is constructed subject to time window $[a_p, b_p']$. |
| (A.12) imposes ride time of each passenger constraints. | Passenger $p$'s network is constructed subject to time window $[a_p, b_p']$. |
| (A.13) imposes capacity constraints. | The structure of state–space–time networks ensures that each vehicle transports a number of passengers up to its capacity at a time. |
| (A.14) defines that the decision variables are binary. | (6) defines binary decision variables. |

## 5 Lagrangian relaxation-based solution approach

Defining multi-dimensional decision variables $y(v,i,j,t,s,w,w')$ leads to computational challenges for the large-scale real-world data sets, which should be addressed properly by specialized programs and an innovative solution framework. We reformulate the problem by relaxing the complicating constraints (5) into the objective function and introducing Lagrangian multipliers, $\lambda(p)$, to construct the dualized Lagrangian function (7).

$$L = \sum_{v \in (V \cup V^*)} \sum_{(i,j,t,s,w,w') \in B_v} c(v,i,j,t,s,w,w') y(v,i,j,t,s,w,w') + \\ \sum_{p \in P} \lambda(p) \left[ \sum_{v \in (V \cup V^*)} \sum_{(i,j,t,s,w,w') \in \Psi_{p,v}} y(v,i,j,t,s,w,w') - 1 \right] \qquad (7)$$

Therefore, the new relaxed problem can be written as follows:

$$Min \; L \tag{8}$$

s.t.

$$\sum_{(i,j,t,s,w,w')\in B_v} y(v,i,j,t,s,w,w') = 1 \quad i = o'_v, t = e_v, w = w' = w_0, \forall v \in (V \cup V^*) \tag{9}$$

$$\sum_{(i,j,t,s,w,w')\in B_v} y(v,i,j,t,s,w,w') = 1 \quad j = d'_v, s = l_v, w = w' = w_0, \forall v \in (V \cup V^*) \tag{10}$$

$$\sum_{(j,s,w'')} y(v,i,j,t,s,w,w'') - \sum_{(j',s',w')} y(v,j',i,s',t,w',w) = 0 \quad (i,t,w) \notin \{(o'_v,e_v,w_0),(d'_v,l_v,w_0)\}, \forall v \in (V \cup V^*) \tag{11}$$

$$y(v,i,j,t,s,w,w') \in \{0,1\} \quad \forall (i,j,t,s,w,w') \in B_v, \forall v \in (V \cup V^*) \tag{12}$$

If we further simplify function $L$, the problem will become a time-dependent least-cost path problem in the constructed state–space–time network. The simplified Lagrangian function $L$ can be written in the following form:

$$L = \sum_{v\in(V\cup V^*)} \sum_{(i,j,t,s,w,w')\in B_v} \xi(v,i,j,t,s,w,w') y(v,i,j,t,s,w,w') - \sum_{p\in P} \lambda(p) \tag{13}$$

Where the generalized arc cost $\xi(v,i,j,t,s,w,w')$ equals $c(v,i,j,t,s,w,w') + \lambda(p)$ for each arc $(i,j,t,s,w,w') \in \Psi_{p,v}$, and equals $c(v,i,j,t,s,w,w')$, otherwise.

### 5.1 Time-dependent forward dynamic programming and computational complexity

Several efficient algorithms have been proposed to compute time-dependent shortest paths in a network with time-dependent arc costs (Ziliaskopoulos and Mahmassani, 1993; Chabini, 1998). In this section, we use a time-dependent dynamic programming (DP) algorithm to solve the least-cost path problem obtained in section 4. The structure of the state–space–time network ensures that time always advances on the arcs of the networks. In this paper, let us consider the unit of time as 1 min. Let $\mathcal{N}$ denote the set of nodes including both physical transportation and dummy nodes, $\mathcal{A}$ denote the set of links, $\mathcal{T}$ denote the set of time stamps covering all vehicles' time horizons, $\mathcal{W}$ denote the set of all feasible passenger carrying states, and $L(i,t,w)$ denote the label of vertex $(i,t,w)$ and term "pred" stands for the predecessor. Algorithm 1 described below uses forward dynamic programming:

```
// Algorithm 1: Time-dependent forward dynamic programming algorithm
    for each vehicle v ∈ (V ∪ V*) do
    begin
        // initialization
        L(.,.,.) := +∞;
        node pred of vertex (.,.,.) := −1;
        time pred of vertex (.,.,.) := −1;
        state pred of vertex (.,.,.) := −1;
        // vehicle v starts its route from the empty state at its origin at the earliest departure time L(o'_v, e_v, w_0) := 0;
        for each time t ∈ [e_v, l_v] do
        begin
            for each link (i, j) do
            begin
                for each state w do
                begin
                    derive downstream state w' based on the possible state transition on link (i, j);
                    derive arrival time s = t + TT(i, j, t);
                    if (L(i, t, w) + ξ(v, i, j, t, s, w, w') < L(j, s, w')) then
                    begin
                        L(j, s, w') := L(i, t, w) + ξ(v, i, j, t, s, w, w') ; //label update
                        node pred of vertex (j, s, w') := i;
                        time pred of vertex (j, s, w') := t;
```

```
                        state pred of vertex (j, s, w') := w;
                    end;
                end;
            end;
        end;
    end;
```

Let's define $|\mathcal{T}|$, $|\mathcal{A}|$, $|\mathcal{W}|$ as the number of time stamps, links, and passenger carrying states, respectively. Therefore, the worst-case complexity of the DP algorithm is $|\mathcal{V}||\mathcal{T}||\mathcal{A}||\mathcal{W}|$, which can be interpreted as the maximum number of steps to be performed in this algorithm in this four-loop structure, corresponding to the sequential loops over vehicle, time, link, and starting carrying state dimensions. It should be remarked that the ending state $w'$ is uniquely determined by the starting state $w$ and the related link $(i,j)$ depending on its service type: pickup, delivery, or pure transportation. In a transportation network, the size of links is much smaller than the counterpart in a complete graph, that is, $|\mathcal{A}| \ll |\mathcal{N}||\mathcal{N}|$; in fact, the typical out-degree of a node in transportation networks is about 2-4.

Table 6 shows detailed comparisons between the existing DP-based approach (Psaraftis, 1983 and Desrosiers et al. 1986) and our proposed method. We guarantee the completeness of state representation. The state representation of Psaraftis (1983), $(L, k_1, k_2, \ldots, k_n)$, consists of $L$, the location currently being visited, and $k_i$, the status of passenger $i$. In this representation, $L = 0$, $L = i$, and $L = i + n$ denote starting location, passenger $i$'s origin, and passenger $i$'s destination, respectively. In addition, the status of passenger $i$ is chosen from the set $\{1,2,3\}$, where 3 means passenger $i$ is still waiting to be picked up, 2 means passenger $i$ has been picked up but the service has not been completed, and 1 means passenger $i$ has been successfully delivered. This cumulative passenger service state representation (in terms of $k_1, k_2, \ldots, k_P$) requires a space complexity of $O(3^P)$, while our proposed (prevailing) passenger carrying state representation has a much smaller space requirement of $\sum_{k=0}^{Cap_v} C_k^P$ when the vehicle capacity is low (e.g. 2 or 3 for taxi). Desrosiers et al. (1986) use state representation $(S, i)$, where $S$ is the set of passengers' origin, $\{1, \ldots, n\}$, and destination, $\{n + 1, \ldots 2n\}$. State $(S, i)$ is defined if and only if there exists a feasible path that passes through all nodes in $S$ and ends at node $i$. In fact, our time-dependent state $(w, i, t)$, which is jointly defined by three indexes: $(i)$ the status of customers, $(ii)$ the current node being visited, and $(iii)$ the current time, is more focused on the time-dependent current state at exact time stamp $t$, while $(L, k_1, k_2, \ldots, k_n)$ and $(S, i)$ representations use a time-lagged time-period-based state representation to cover complete or mutually exclusive states from time 0 to time $t$.

**Table 6**
Comparison between existing DP based approach and the method proposed in this paper.

| Features | Existing DP based approach | | DP proposed in this paper |
|---|---|---|---|
| | Psaraftis (1983) | Desrosiers et al. (1986) | |
| Type of problem | Single vehicle, Many-to-many, Single depot | Single vehicle, Many-to-many, Single depot | Multiple vehicle, Many-to-many, Multiple depot |
| Network | Consists of passengers' origin and destination nodes and the vehicle depot | Consists of passengers' origin and destination nodes and the vehicle depot | Consists of transportation nodes, passengers' origin and destination, and vehicles' depots |
| Time-dependent link travel time | No | No | Yes |
| Objective function | Minimize route duration | Minimize total distance traveled | Minimize total routing cost consisting of transportation and waiting costs |
| State | State–space $(L, k_1, k_2, \ldots, k_n)$ | State–space $(S, i)$ | State–space–time $(w, i, t)$ |
| Stage | Node index | Node index | Time index |
| States reduction due to the vehicle capacity and time windows | Yes | Yes | Yes |

We come back to the illustrative example presented in section 3.1. Let's assume the routing cost of a transportation or service arc traversed by a physical vehicle is $22/h, while the routing cost of a transportation or service arc traversed by a virtual vehicle is $50/h. Moreover, assume that the waiting cost of a physical vehicle is $15/h, while the waiting cost of a virtual vehicle is assumed to be $0/h. Table 7 shows how the label of each vertex is calculated by the DP solution algorithm presented above. Note that $w_0$, $w_1$, $w_2$, and $w_3$ are passenger carrying states [ _ _ ], [ $p_1$ _ ], [ $p_1$ $p_2$ ], and [ _ $p_2$ ], respectively. For instance, according to Fig. 1, traveling from node 4 to node 2 takes 2 min. Since the number written on each link denotes the time-invariant travel time $TT(i,j)$, we can conclude that travel time for link (4,2) starting at time stamp $t = 2$ is also 2 min. To update the label corresponding to node 2, it is sufficient to calculate the routing cost of the stated arc in terms of dollars which can be obtained by $\frac{2}{60} \times 22(\frac{\$}{h}) = 0.73(\$)$ and add it to the current label of node 4 which is 0.37. Therefore, the updated label for node 2 will be 1.1. Similarly, we can calculate the routing cost of a waiting link $(o_2, o_2)$ starting at time stamp $t = 7$ by $\frac{1}{60} \times 15\left(\frac{\$}{h}\right) = 0.25$ ($).

**Table 7**
State–space–time trajectory for ride-sharing service trip with node sequence $(o'_1, 4, 2, o_1, 2, o_2, 2, 5, 6, 3, d_1, 3, d_2, 3, 1, d'_1)$.

| Time index | 1 | 2 | 4 | 5 | 6 | 7 | 8 | 9 | 10 | 11 | 12 | 13 | 14 | 15 | 16 | 18 | 19 | 20 |
|---|---|---|---|---|---|---|---|---|---|---|---|---|---|---|---|---|---|---|
| Node index | $o'_1$ | 4 | 2 | $o_1$ | 2 | $o_2$ | $o_2$ | 2 | 5 | 6 | 3 | $d_1$ | 3 | $d_2$ | 3 | 1 | $d'_1$ | $d'_1$ |
| State index | $w_0$ | $w_0$ | $w_0$ | $w_1$ | $w_1$ | $w_2$ | $w_2$ | $w_2$ | $w_2$ | $w_2$ | $w_2$ | $w_2$ | $w_3$ | $w_3$ | $w_0$ | $w_0$ | $w_0$ | $w_0$ |
| Cost | 0.0 | .37 | .73 | .37 | .37 | .37 | .25 | .37 | .37 | .37 | .37 | .37 | .37 | .37 | .37 | .37 | .73 | .37 | 0.0 |
| Cumulative cost | 0.0 | .37 | 1.1 | 1.47 | 1.84 | 2.21 | 2.46 | 2.83 | 3.2 | 3.57 | 3.94 | 4.31 | 4.68 | 5.05 | 5.42 | 6.15 | 6.52 | 6.52 |

### 5.2 Lagrangian relaxation-based solution procedure

In this section, we describe the Lagrangian relaxation (LR) solution approach implemented to solve the time-dependent least cost path problem presented in section 5. According to Eq. (13), $\xi(v, i, j, t, s, w, w')$ is only updated for $\forall (v, i, j, t, s, w, w') \in \Psi_{p,v}$. Table 8 lists the notations for the sets, indices and parameters required for the Lagrangian relaxation algorithm.

**Table 8**
Notations used in LR algorithm.

| Symbol | Definition |
|---|---|
| $\lambda^k(p)$ | Lagrangian relaxation multiplier corresponding to the passenger $p$'s pickup request constraint at iteration $k$ |
| $\xi(v, i, j, t, s, w, w')$ | Modified routing cost of arc $(i, j, t, s, w, w')$ after introducing Lagrangian multipliers |
| $k$ | Iteration number |
| $Y$ | Set of solution vectors $y(v, i, j, t, s, w, w')$ |
| $LB^k$ | Global lower bound for the primal problem at iteration $k$ |
| $UB^k$ | Global upper bound for the primal problem at iteration $k$ |
| $Y_{LB}^k$ | Set of lower bound solution vectors $Y$ at LR iteration $k$ |
| $Y_{UB}^k$ | Set of upper bound solution vectors $Y$ at LR iteration $k$ |
| $\theta^k$ | Step size at iteration $k$ |
| $LB^*$ | Best global lower bound for the primal problem |
| $UB^*$ | Best global upper bound for the primal problem |
| $Y^*$ | Best solution vectors derived from the best lower bound |
| $base\_profit$ | The amount of money (in terms of dollars) passenger $p$ initially offers to be served |

The optimal value of the Lagrangian dual problem provides a lower bound for the primal problem. To find the optimal solution for the Lagrangian dual problem, it is sufficient to compute time-dependent least cost state–space–time path for each vehicle $v$ based on updated arc cost $\xi$'s by calling time-dependent forward dynamic programming algorithm mentioned before.

The optimal solution of the Lagrangian dual problem may or may not be feasible for the primal problem. If the optimal solution of the Lagrangian dual problem is feasible for the primal problem, we have definitely obtained the optimal solution of the primal problem. If not, we apply a heuristic to find an upper bound for the primal solution. In this heuristic, the physical vehicles initially leave their depots to serve unserved customers provided that the money obtained in return for services overweigh the cost of transportation. Finally, if there is any unserved customer remained in the system, in order to avoid infeasibility, the virtual vehicle corresponding to the unserved customer departs from its depot to serve the passenger. The Lagrangian relaxation algorithm can be described as follows:

// Algorithm 2: Lagrangian relaxation algorithm
    // step 0. initialization
- set iteration $k = 0$;
- initialize $Y_{LB}^0, Y_{UB}^0, Y^*$, and $\lambda^0(p)$ to zero;
- initialize $\theta^0(p)$ to $base\_profit$;
- initialize $LB^*$ to $-\infty$; and $UB^*$ to $+\infty$;
- define a termination condition such as if $k$ becomes greater than a predetermined maximum iteration number, or if the relative gap percentage between $LB^*$ and $UB^*$ becomes less than a predefined gap (i.e. 5%);

while termination condition is false, for each LR iteration $k$ do
begin
- reset the visit count for each arc $(v, i, j, t, s, w, w') \in \Psi_{p,v}$ to zero; // $v \in (V \cup V^*)$

// step 1. generating $LB^k$
    // step 1.1. least cost path calculation for each vehicle sub-problem
- initialize $LB^k$ to 0;

for each vehicle $v \in (V \cup V^*)$ do
begin
    // input: $\xi(v, i, j, t, s, w, w')$
- compute time-dependent least cost state–space–time path for vehicle $v$ based on updated arc cost $\xi$'s by calling Algorithm 1;
- update the visit count for each arc $(v, i, j, t, s, w, w') \in \Psi_{p,v}$;
    // output: $Y_{LB}^k$
end;

// step 1.2. update $LB^*$
- update $LB^k$ by substituting solution vector $Y_{LB}^k$ in the objective function of the dual problem (Eq. (13));
- update $LB^*$ by $max(LB^k, current\ LB^*)$ and $Y^*$ by its corresponding solution;

// step 1.3. sub-gradient calculation
- calculate the total number of visits of passenger $p$'s origin by expression (14);

$$\sum_{v \in (V \cup V^*)} \sum_{(i,j,t,s,w,w') \in \Psi_{p,v}} y(v, i, j, t, s, w, w') \tag{14}$$

- compute sub-gradients by Eq. (15);

$$\nabla L_{\lambda^k(p)} = \sum_{v \in (V \cup V^*)} \sum_{(v,i,j,t,s,w,w') \in \Psi_{p,v}} y(v, i, j, t, s, w, w') - 1 \text{ for } \forall p \tag{15}$$

- update arc multipliers by Eq. (16);

$$\lambda^{k+1}(p) = \lambda^k(p) + \theta^k(p) \nabla L_{\lambda^k(p)} \text{ for } \forall p \tag{16}$$

- update arc cost $\xi(v, i, j, t, s, w, w')$ for each arc $(v, i, j, t, s, w, w') \in \Psi_{p,v}$ by Eq. (17);

$$\xi(v, i, j, t, s, w, w') = c(v, i, j, t, s, w, w') + \lambda^{k+1}(p) \tag{17}$$

- update step size by Eq. (18);

$$\theta^{k+1}(p) = \frac{\theta^0(p)}{k+1} \tag{18}$$

// Step 2. generating $UB^k$
    // step 2.1. finding a feasible solution for the primal problem
- set $UB^k = 0$;
- adopt the passenger-to-vehicle assignment matrix from the lower bound solution in step 1.2;

for each vehicle $v \in (V \cup V^*)$ do
begin
    // if passenger $p$ is served by multiple vehicles, then designate one of the vehicles (e.g. first in the set) to serve this passenger, which means that the other vehicles should not serve this passenger in the upper bound generation stage.
if (passenger $p$ is assigned to physical vehicle $v$) do
begin

        if passenger $p$ has not been already served by any other vehicle
            set arc cost on the pickup arc for passenger $p$ temporarily to $-M$;
            // $M$ is chosen a big positive number in order to attract vehicle $v$ for serving passenger $p$
        else
            set arc cost on the pickup arc for passenger $p$ temporarily to $+M$;
            // $M$ is chosen a big positive number in order to guarantee vehicle $v$ does not serve passenger $p$
        end;
        // if passenger $p$ is not served by any physical vehicle, then designate the corresponding virtual vehicle to serve this passenger.
        if (passenger $p$ is not served by any physical vehicle & vehicle $v$ is the corresponding virtual vehicle for passenger $p$)
            set arc cost on the pickup arc for passenger $p$ temporarily to $-M$;
            // $M$ is chosen a big positive number in order to attract vehicle $v$ for serving passenger $p$
- compute time-dependent least cost path for vehicle $v$ by calling Algorithm 1;
- compute the actual transportation costs (denoted as $TC_v$) along the path solution for vehicle $v$
- update upper bound objective function as $UB^k = UB^k + TC_v$.

    end;
    // The result of this passenger-to-vehicle assignment updating is that each passenger is served by exactly one vehicle (either physical or virtual).
// step 2.2. update $UB^k$
- update $UB^k$ by substituting solution vector $Y_{UB}^k$ in the objective function of the primal problem;
// step 2.3. update $UB^*$
- $UB^* = min(UB^k, current\ UB^*)$;
- find the relative gap percentage between $LB^*$ and $UB^*$ by $\frac{UB^* - LB^*}{UB^*} \times 100$;
- $k = k + 1$;

end;

We would like to make remarks in following two cases:

    ($i$) In the upper bound solution, all passengers are only served by the physical vehicles. In this case, we can be sure that the total number of physical vehicles has been sufficient to serve all requests. Accordingly, the service prices in the corresponding lower bound solution typically have been set such that the money obtained in return for services overweighs the cost of transportation so that physical vehicles are dispatched to serve customers.

    ($ii$) In the final optimal solution, there might be some passengers who are served by virtual vehicles. Obviously, serving a passenger by a virtual vehicle is expensive due to its transportation cost. In addition, when the virtual vehicle drops off the passenger, it should perform a deadheading trip with significantly high cost from the passenger's destination to its depot (the passenger's origin).

### 5.3 Search region reduction

In this section, we describe how to reduce the search region by the aid of some simple heuristics in which some rational rules are applied. Let $EDT$, $LDT$, $EAT$, and $LAT$ denote the earliest departure time from origin, latest departure time from origin, earliest arrival time to destination, and latest arrival time to destination, respectively. In addition, let $TTSP_{x \to y}$ denote the travel time corresponding to the shortest path from node $x$ to node $y$.

Rule 1. No overlapping time windows: The first rational rule is that if $LAT(p_1) < EDT(p_2)$, then passenger $p_1$ and $p_2$'s ride-sharing is impossible. Fig. 7 illustrates an example of two passengers whose ride-sharing is impossible due to no overlapping time windows.

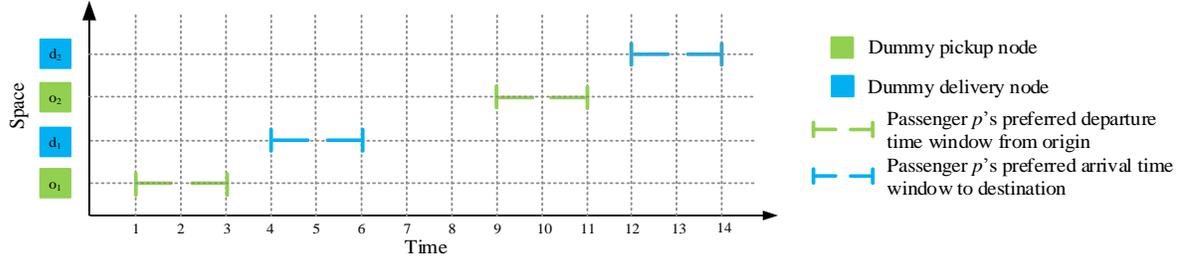

**Fig. 7.** Illustration of the first rational rule for search region reduction.

Rule 2. Travel time is insufficient: The second rational rule can be stated as follows: if $\{LDT(p_2) - EDT(p_1) < TTSP_{o_{p_1} \to o_{p_2}}$ & $LDT(p_1) - EDT(p_2) < TTSP_{o_{p_2} \to o_{p_1}}\}$, then passenger $p_1$ and $p_2$ cannot share their ride with each other. It means that if the maximum time a vehicle can have to go from passenger $p_1$'s origin to $p_2$'s origin, $LDT(p_2) - EDT(p_1)$, is less than the total travel time corresponding to the shortest path from $o_{p_1}$ to $o_{p_2}$, and also if the maximum time a vehicle can have to go from passenger $p_2$'s origin to $p_1$'s origin, $LDT(p_1) - EDT(p_2)$, is less than the total travel time corresponding to the shortest path from $o_{p_2}$ to $o_{p_1}$, then passenger $p_1$ and $p_2$'s ride-sharing is impossible. Similarly, if $\{LAT(p_2) - EAT(p_1) < TTSP_{d_{p_1} \to d_{p_2}}$ & $LAT(p_1) - EAT(p_2) < TTSP_{d_{p_2} \to d_{p_1}}\}$, then passenger $p_1$ and $p_2$'s ride-sharing is impossible. The total number of passenger carrying states is dramatically decreased via this rule. Fig. 8 illustrates the second rule by an example. Suppose two requests with two origin–destination pairs should be served by a vehicle. Fig. 8(a) illustrates transportation network with the corresponding dummy nodes and time windows. According to the Fig. 8(a), $TTSP_{o_{p_1} \to o_{p_2}}$ and $TTSP_{o_{p_2} \to o_{p_1}}$ are 5 and 6, respectively. Since $\{(6-4) < 5$ & $(5-4) < 6\}$, then passenger $p_1$ and $p_2$'s ride-sharing is impossible.

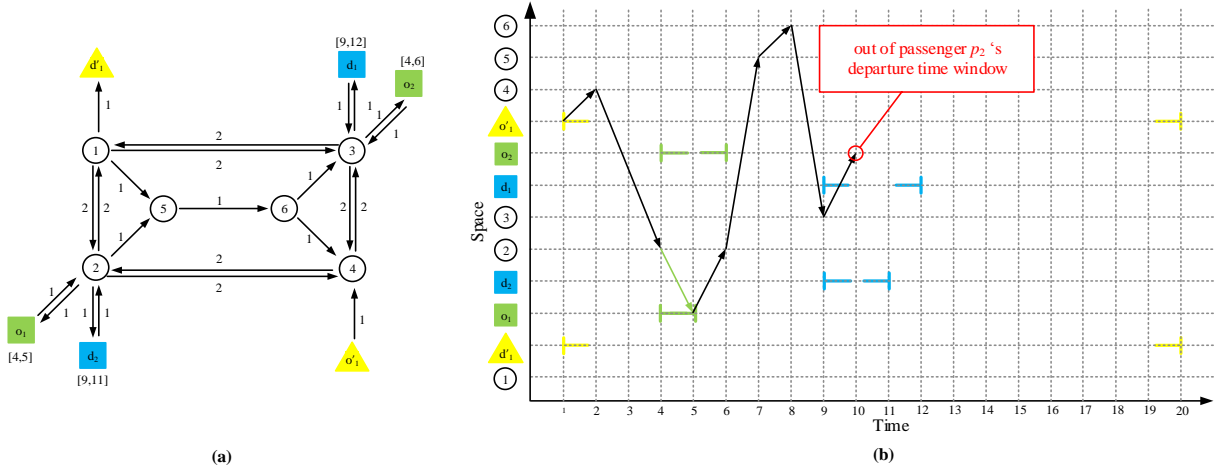

**Fig. 8.** Illustration of the second rational rule for search region reduction; (a) transportation network with the corresponding dummy nodes and time windows; (b) vehicle 1's space–time network.

Rule 3. A node is too far away from the vehicle starting or ending depot: The third rational rule is stated as follows: if $(TTSP_{O_v \to x} + TTSP_{x \to d_v}) > (LAT(v) - EDT(v))$, then vehicle $v$ does not have enough time to visit node $x$ in its time horizon; therefore, node $x$ is not accessible for vehicle $v$ and should not be considered in vehicle $v$'s search region. Note that node $x$ can be any physical or dummy node. Fig. 9 illustrates the third rule by an example. Suppose a passenger with an origin–destination pair should be served by a vehicle. Fig. 9(a) illustrates transportation network with the corresponding dummy nodes and time windows. Fig. 9(b) shows that passenger $p_1$'s origin, $o_1$, is not accessible for the vehicle. In addition to this rule, we can also say that a passenger is inaccessible for a vehicle if the time for a vehicle to pick up the passenger and visit his destination is longer than the vehicle's time window.

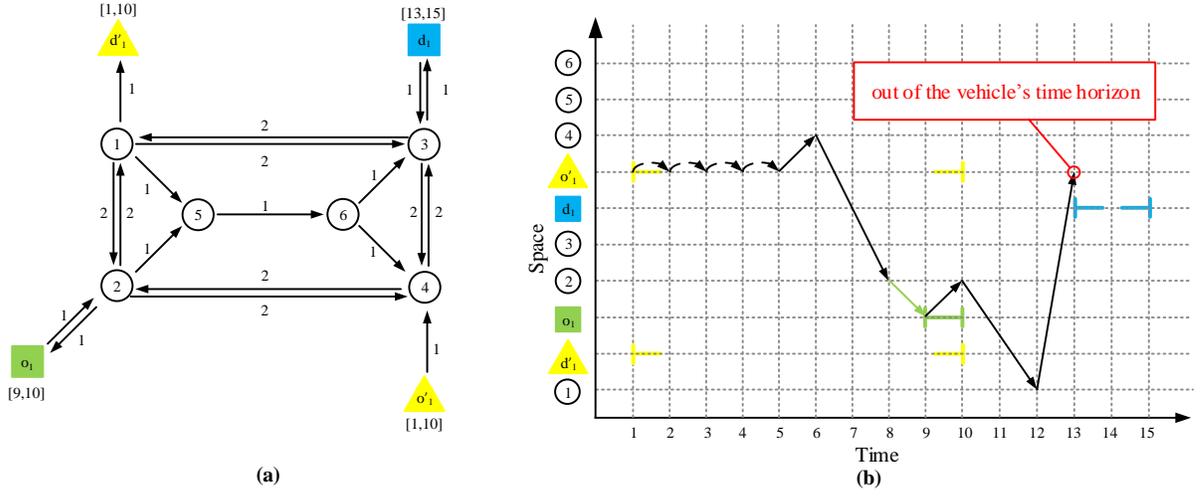

**Fig. 9.** Illustration of the third rational rule for search region reduction; (a) transportation network with the corresponding dummy nodes and time windows; (b) vehicle 1's space–time network.

The first three rules are hard rules at which we are able to eliminate some vertexes in the state–space–time networks. The forth heuristic is the way of estimating the search region reduction ratio. Let path $\alpha$ be the longest possible path in vehicle $v$'s state–space–time networks with total travel time $\tau_\alpha$. Let $m_p$ denote the middle point of passenger $p$'s departure time window. Therefore, $m_p = \frac{EDT(p)+LDT(p)}{2}$. Let's assume that $M$, the middle point of a passenger's departure time window, is a random variable uniformly distributed in vehicle $v$'s time horizon with $LAT(v) - EDT(v)$ length. It may be reasonable to assume that if $|m_{p_1} - m_{p_2}| > \tau_\alpha$, then passenger $p_1$ and $p_2$ cannot be in the same vehicle at a time. We use an example to show that this rule can reduce the search region considerably. Assume vehicle $v$'s time window is [0, 240], and $M$ is a random variable uniformly distributed in vehicle $v$'s time horizon [0,240]. Let's assume $\tau_\alpha = 60$ min. The probability of having two passengers who share their ride with each other can be calculated by finding the $Prob(|m_{p_1} - m_{p_2}| \leq 60\ minutes)$, where $m_{p_1}$ and $m_{p_2}$ are randomly generated from [0, 240]. This probability equals to $\frac{7}{16} = 43.75\%$. This can be shown with the following derivation. The shaded area in in Fig. 10 shows $Prob(|m_{p_1} - m_{p_2}| \leq 60\ )$.

$$Prob(|m_{p_1} - m_{p_2}| \leq 60\ ) = Prob(-60 \leq m_{p_1} - m_{p_2} \leq 60\ )$$
$$Prob(-60 \leq m_{p_1} - m_{p_2} \leq 60\ ) = 1 - [Prob(m_{p_1} - m_{p_2} < -60\ ) + Prob(m_{p_1} - m_{p_2} > 60\ )]$$
$$= 1 - \left[\frac{\frac{180 \times 180}{2}}{240 \times 240} + \frac{\frac{180 \times 180}{2}}{240 \times 240}\right] = \frac{7}{16}$$

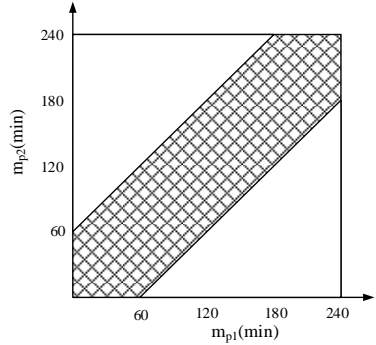

**Fig. 10.** The probability of having two passengers who share their ride with each other where $m_{p_1}$ and $m_{p_2}$ are uniformly distributed in [0, 240]. Note that $\tau_\alpha = 60$ min.

Therefore, by considering this practical rule in this example, we can reduce the total number of passenger carrying states in which two passengers share their ride with each other by more than half. By considering this rational rule, calculating the probability of having more than two passengers at the same time in vehicle $v$ is more complicated, but at least we know that the probability of having $k$ number of passengers ($k > 2$) who may share their ride with each other is certainly less than 43.75%.

## 6 Computational results and discussions

The algorithms described in this paper were coded in C++ platforms. The experiments were performed on an Intel Workstation running two Xeon E5-2680 processors clocked at 2.80 GHz with 20 cores and 192 GB RAM running Windows Server 2008 x64 Edition. In addition, parallel computing and OpenMP technique (Chandra et al., 2000) are implemented for generating lower bound and upper bound at each iteration in the Lagrangian relaxation algorithm. In this section, we initially examine our proposed model on a six-node transportation network followed by the medium-scale and large-scale transportation networks, Chicago and Phoenix, to demonstrate the computational efficiency and solution optimality of our developed algorithm. The scenarios and test cases are randomly generated in those transportation networks. Moreover, we test our algorithms on the modified version of instances proposed by Ropke and Cordeau (2009) which is publicly available at http://www.diku.dk/~sropke/.

### 6.1 Illustrative cases

As we mentioned in section 5.1, it is assumed that the routing cost of a transportation or service arc traversed by a physical vehicle is $22/h, while the routing cost of a transportation or service arc traversed by a virtual vehicle is $50/h. Moreover, the waiting cost of a physical vehicle either at a transportation or at a dummy node is $15/h (waiting at dummy nodes corresponding starting and ending depots has $0/h cost), while the waiting cost of a virtual vehicle at any node is assumed to be $0/h. The value of $base\_profit$ is also assumed to be $10 for all passengers. Initially, we test our algorithm on the six-node transportation network illustrated in Fig. 1(a) for six scenarios. Table 9 shows these scenarios with various number of passengers and vehicles, origin–destination pairs, and passengers' departure and arrival time windows. Then, we will examine the results corresponding to each scenario individually. Terms "TW" and "TH" stands for time window and time horizon, respectively. Table 10 shows the results corresponding each scenario.

Scenario I. Two passengers are served by one vehicle, where passengers have different origin–destination pairs with overlapping time windows. In this case, the vehicle serves both passengers in their preferred time windows through ride-sharing mode.

Scenario II. Two passengers with different origin–destination pairs are served by one vehicle; however, unlike in scenario I, passengers could not share their ride with each other due to their time windows. In this case, the vehicle may wait at any node to finally serve both passengers.

Scenario III. Two passengers with different origin–destination pairs and one vehicle are present in the system; however, due to the passengers' overlapping time windows, serving both passengers by one vehicle is impossible. Therefore, the driver would prefer to transport a passenger incurring the least cost. In this case, passenger $p_1$ is selected to be served.

Scenario IV. Two passengers with different origin–destination pairs and two vehicles are present in the system and, due to the passengers' and vehicles' time windows, $p_1$ is assigned to $v_1$ and $p_2$ is assigned to $v_2$.

Scenario V. Three passengers are served by one vehicle, where passengers have different origin–destination pairs with overlapping time windows. In this case, the vehicle serves all passengers in their preferred time windows through ride-sharing mode.

Scenario VI. One passenger and two vehicles are present in the system. In this case, two vehicles compete for serving the passenger. Ultimately, the vehicle whose routing is less costly wins the competition and serves the passenger.

**Table 9**
Six scenarios with various numbers of passengers and vehicles, origin–destination pairs, and passengers' departure and arrival time windows.

| Scenario | I | II | III | IV | V | VI |
|---|---|---|---|---|---|---|
| Number of passengers | 2 | 2 | 2 | 2 | 3 | 1 |
| Number of vehicles | 1 | 1 | 1 | 2 | 1 | 2 |
| $o_1$ | Node 2 | Node 2 | Node 2 | Node 2 | Node 2 | Node 2 |
| $d_1$ | Node 6 | Node 6 | Node 1 | Node 1 | Node 3 | Node 6 |
| $o_2$ | Node 5 | Node 5 | Node 3 | Node 3 | Node 5 | – |
| $d_2$ | Node 3 | Node 3 | Node 6 | Node 6 | Node 3 | – |
| $o_3$ | – | – | – | – | Node 6 | – |
| $d_3$ | – | – | – | – | Node 1 | – |
| $o'_1$ | Node 4 | Node 4 | Node 4 | Node 2 | Node 4 | Node 4 |
| $d'_1$ | Node 1 | Node 1 | Node 1 | Node 1 | Node 1 | Node 1 |
| $o'_2$ | – | – | – | Node 3 | – | Node 6 |
| $d'_2$ | – | – | – | Node 6 | – | Node 1 |
| $TW_{o_1}$ | [5, 7] | [5, 7] | [4, 5] | [4, 5] | [4, 7] | [4, 7] |
| $TW_{d_1}$ | [9, 12] | [9, 12] | [8, 10] | [8, 10] | [13, 16] | [9, 12] |
| $TW_{o_2}$ | [8, 10] | [16, 19] | [3, 5] | [4, 6] | [7, 10] | – |
| $TW_{d_2}$ | [11, 14] | [21, 24] | [11, 14] | [11, 14] | [14, 18] | – |
| $TW_{o_3}$ | – | – | – | – | [10, 13] | – |
| $TW_{d_3}$ | – | – | – | – | [19, 23] | – |
| $TH_{v_1}$ | [1, 30] | [1, 30] | [1, 30] | [1, 30] | [1, 30] | [1, 30] |
| $TH_{v_2}$ | – | – | – | [1, 30] | – | [1, 30] |

**Table 10**
Results obtained from testing our algorithm on the six-node transportation network for six scenarios.

| iteration $k$ | $LB^*$ | $UB^*$ | gap% | vehicles assigned to $p_1, p_2$, and $p_3$ | $\lambda^k(p_1)$ | $\lambda^k(p_2)$ | $\lambda^k(p_3)$ |
|---|---|---|---|---|---|---|---|
| Scenario I. Two passengers are served by one vehicle through ride-sharing mode. | | | | | | | |
| 1 | 1.47 | 5.75 | 74.5% | $v_1, v_1, -$ | 10 | 10 | – |
| 2 | 1.47 | 5.75 | 74.5% | $v_1, v_1, -$ | 5 | 5 | – |
| 3 | 5.75 | 5.75 | 0.0% | $v_1, v_1, -$ | 5 | 5 | – |
| Scenario II. Two passengers are served by one vehicle (not through ride–sharing mode). | | | | | | | |
| 1 | 1.47 | 7.22 | 79.68% | $v_1, v_1, -$ | 10 | 10 | – |
| 2 | 5.55 | 7.22 | 23.10% | $v_1, v_1, -$ | 5 | 5 | – |
| 3 | 7.22 | 7.22 | 0.0% | $v_1, v_1, -$ | 5 | 5 | – |
| Scenario III. Two passengers and one vehicle; one passenger remains unserved. | | | | | | | |
| 1 | 1.47 | 10.43 | 85.94% | $v_1, v_2^*, -$ | 10 | 10 | – |
| 2 | 7.1 | 10.43 | 31.95% | $v_1, v_2^*, -$ | 5 | 10 | – |
| 3 | 10.43 | 10.43 | 0.0% | $v_1, v_2^*, -$ | 5 | 10 | – |
| Scenario IV. Two passengers and two vehicles; each vehicle is assigned to a passenger | | | | | | | |
| 1 | 2.2 | 6.13 | 64.13% | $v_1, v_2, -$ | 10 | 10 | – |
| 2 | 2.2 | 6.13 | 64.13% | $v_1, v_2, -$ | 5 | 5 | – |
| 3 | 6.13 | 6.13 | 0.0% | $v_1, v_2, -$ | 5 | 5 | – |
| Scenario V. Three passengers are served by one vehicle through ride-sharing mode | | | | | | | |
| 1 | 1.47 | 6.97 | 78.95% | $v_1, v_1, v_1$ | 10 | 10 | 10 |
| 2 | 1.47 | 6.97 | 78.95% | $v_1, v_1, v_1$ | 5 | 5 | 5 |
| 3 | 6.97 | 6.97 | 0.0% | $v_1, v_1, v_1$ | 5 | 5 | 5 |
| Scenario VI. Two vehicles compete for serving a passenger | | | | | | | |
| 1 | 2.57 | 5.13 | 50.0% | $v_1, -, -$ | 10 | – | – |
| 2 | 2.63 | 5.13 | 48.70% | $v_1, -, -$ | 10 | – | – |
| 3 | 5.13 | 5.13 | 0.0% | $v_1, -, -$ | 10 | – | – |

Fig. 11 also presents the vehicle routing corresponding each scenario. We increase the number of passengers and vehicles to show the computational efficiency and solution optimality of our developed algorithm. Table 11 shows the results for the six-node transportation network when the numbers of passengers and vehicles have been increased.

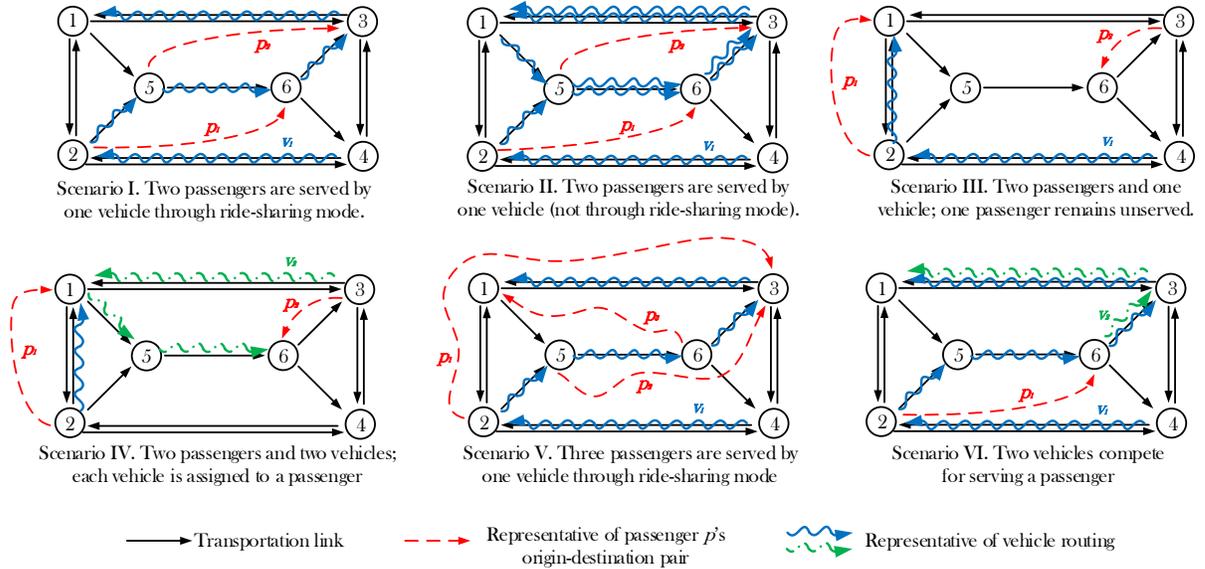

Fig. 11. The vehicle routing corresponding each scenario.

**Table 11**
Results for the six-node transportation network.

| Test case number | Number of iterations | Number of passengers | Number of vehicles | $LB^*$ | $UB^*$ | Gap (%) | Number of passengers not served | CPU running time (s) |
|---|---|---|---|---|---|---|---|---|
| 1 | 30 | 6 | 1 | 15.83 | 15.83 | 0.00% | 0 | 5.94 |
| 2 | 30 | 12 | 2 | 33.17 | 33.17 | 0.00% | 0 | 12.02 |
| 3 | 30 | 24 | 4 | 61.67 | 65.33 | 5.61% | 0 | 30.97 |

We explain the pricing mechanism in this algorithm via test case 1 with 6 passengers and 1 vehicle. Fig. 12 shows $\lambda^k(p_i)$, $i = 1, 2, \ldots, 6$, along 30 iterations. It is clear that each passenger's Lagrangian multiplier ultimately converges to a specific value. This value can be literally interpreted as the passenger $p$'s service price. Through the pricing mechanism of this algorithm, the provider would be able to offer a reasonable bid to its customers to be served.

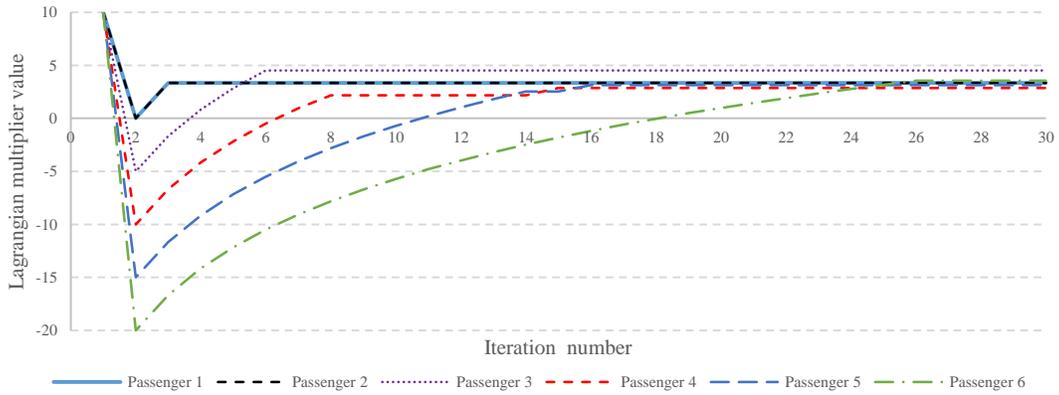

**Fig. 12.** Lagrangian multipliers along 30 iterations in test case 1 for the six-node transportation network.

### 6.2 Results from medium-scale and large-scale transportation networks

In our computational experiments for the medium-scale and large-scale networks, for simplicity, we assume that each passenger has a fixed departure time (the earliest and latest departure time are the same). In addition, we assume that no passenger has a preferred time window for arrival to his destination. Tables 12 and 13 show the results for the Chicago transportation network, shown as Fig. 13(a) with 933 nodes and 2,967 links, and the Phoenix transportation network, as shown in Fig. 13(b) with 13,777 nodes and 33,879 links, respectively.

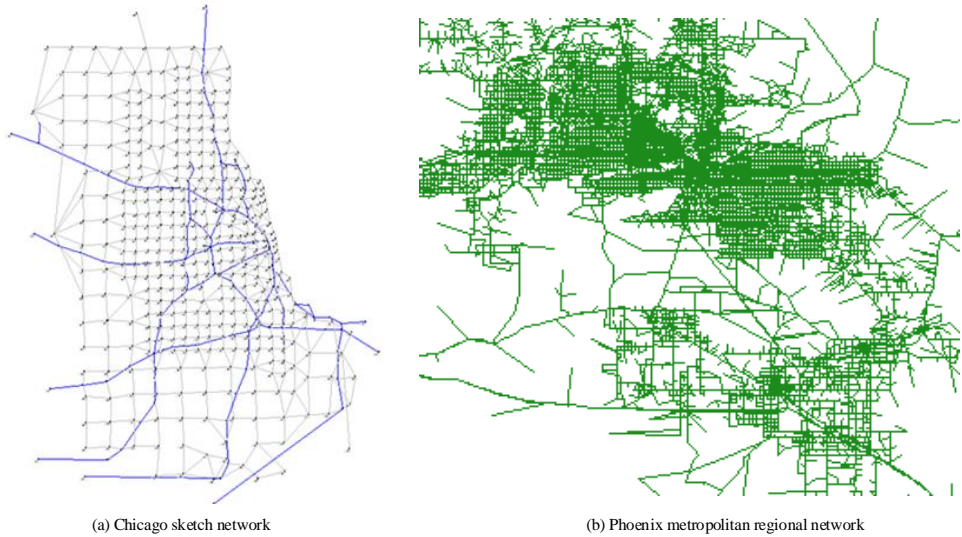

(a) Chicago sketch network        (b) Phoenix metropolitan regional network

**Fig. 13.** Medium and large-scale transportation networks for computational performance testing.

Note that we generally run the algorithm for a fixed number of iterations; however, the algorithm may converge in less number of iterations. Fig. 14 shows the gap percentage along 20 iterations corresponding each test case.

**Table 12**
Results for the Chicago network with 933 transportation nodes and 2,967 links.

| Test case number | Number of iterations | Number of passengers | Number of vehicles | $LB^*$ | $UB^*$ | Gap (%) | Number of passengers not served | CPU running time (s) |
| --- | --- | --- | --- | --- | --- | --- | --- | --- |
| 1 | 20 | 2 | 2 | 108.43 | 108.43 | 0.00% | 0 | 17.43 |
| 2 | 20 | 11 | 3 | 352.97 | 352.97 | 0.00% | 0 | 91.87 |
| 3 | 20 | 20 | 5 | 616.66 | 626.18 | 1.52% | 1 | 327.51 |
| 4 | 20 | 46 | 15 | 1586.81 | 1664.07 | 4.64% | 2 | 4681.52 |
| 5 | 20 | 60 | 15 | 1849.98 | 1878.55 | 1.52% | 3 | 7096.50 |

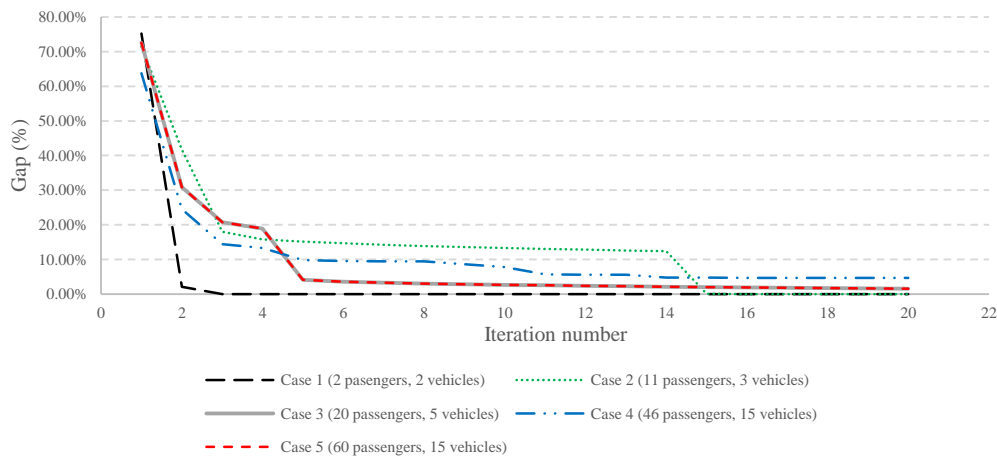

**Fig. 14.** Gap percentage along 20 iterations corresponding each test case in Chicago network.

As you can see in Fig. 14, after 10–15 iterations, the sub-gradient algorithm is typically able to converge to a small gap (about 5%) for the Chicago Network.

**Table 13**
Results for the Phoenix network with 13,777 transportation nodes and 33,879 links.

| Test case number | Number of iterations | Number of passengers | Number of vehicles | $LB^*$ | $UB^*$ | Gap (%) | Number of passengers not served | CPU running time (s) |
|---|---|---|---|---|---|---|---|---|
| 1 | 6 | 4 | 2 | 70.95 | 70.95 | 0.00% | 0 | 110.39 |
| 2 | 6 | 10 | 5 | 191.55 | 207.05 | 7.49% | 1 | 398.37 |
| 3 | 6 | 20 | 6 | 310.37 | 310.37 | 0.00% | 0 | 1323.18 |
| 4 | 6 | 40 | 12 | 622.23 | 622.23 | 0.00% | 0 | 3756.51 |
| 5 | 6 | 50 | 15 | 784.07 | 784.07 | 0.00% | 0 | 6983.19 |

### 6.3 Handling randomly generated test instances

To further examine the computational efficiency and solution optimality of our proposed algorithm, we also test our algorithms on randomly generated instances proposed by Ropke and Cordeau (2009) which is publicly available at http://www.diku.dk/~sropke/. The data set introduced by Ropke and Cordeau (2009) is the modified version of instances employed by Ropke et al. (2007) initially proposed by Savelsbergh and Sol (1998). In this data set, the coordinates of passengers' pickup and drop-off locations are randomly selected and uniformly distributed over a $[0,50] \times [0,50]$ square. In addition, they considered a single depot located in the center of the square. The load $q_i$ of passenger $i$ is randomly selected from $[5, Q]$, where $Q$ is the maximum capacity of the vehicle. A planning horizon $[0,600]$ is considered. Feasible departure and arrival time windows are also randomly generated for each passenger.

Ropke and Cordeau (2009) formulate the PDPTW on a network that is built based on demand request nodes, and the links are defined as direct connections between pickup and delivery nodes (without explicitly considering transportation links or paths), while we formulate the PDPTW on transportation networks. To test our algorithms on their data set, we need to convert their origin–destination network to a transportation network. Specifically, we treat their demand node-oriented network as a transportation network, and each origin/destination node acts as a transportation node. As a result, in this converted transportation network, each transportation node is connected to all other transportation nodes, and similar to what we performed before, dummy nodes are added and connected to their corresponding transportation nodes. Obviously, the constructed transportation network is a complete digraph with a very large number of links. In this data set, the coordinates of passengers' pickup and drop-off locations are randomly chosen and uniformly distributed over a small square, so the densely distributed passengers could impose a difficult problem of assigning different vehicles to different passengers. In comparison, in our test data set, the Chicago and Phoenix transportation networks, the vehicles and passengers are naturally spatially and sparsely distributed such that fewer vehicles compete for serving a particular passenger. Thus, in that situation, the first stage of the vehicle assignment problem could be easily solved using our proposed Lagrangian relaxation framework with a good matching between the vehicles and passengers.

**6.4 Challenges of multi-vehicle assignment problems and usefulness of single-vehicle routing algorithm**

In general, VRPPDTW even for the single vehicle cases is still categorized as one of the toughest tasks of combinatorial optimization (Azi et al., 2007; Hernández-Pérez and Salazar-González, 2009; and Häme, 2011). Several approaches have been recently suggested to resolve the issue mentioned above by converting multi-vehicle cases to the single-vehicle ones. For instance, Häme and Hakula (2015) have suggested a maximum cluster algorithm in which the multi-vehicle solution is based on a recursive single-vehicle algorithm.

To fully address the complexity of assigning different vehicles to multiple passengers, Fisher et al. (1997) proposed a new set of variables $x(p,v)$ and decomposed constraints (5) to two sets of constraints: constraints (19) and (20).

$$\sum_{v \in (V \cup V^*)} x(v,p) = 1 \qquad \forall p \quad (19)$$
$$\sum_{(j,s,w')} y(v,i,j,t,s,w,w') = x(v,p) \qquad (i,j,t,s,w,w') \in \Psi_{p,v}, \forall p, \forall v \quad (20)$$

Constraints (19) guarantee that each passenger is visited exactly once. Constraints (20) control vehicle $v$'s route and show the relations between the variables $x(v,p)$ and $y(v,i,j,t,s,w,w')$. By using the Lagrangian decomposition method and relaxing these two sets of constraints into the objective function, the main problem can be decomposed to two sub-problems where sub-problem (1) becomes a semi-assignment problem and sub-problem (2) as a time-dependent least-cost path problem. The first sub-problem can be easily solved by inspection. The second sub-problem also can be solved by computationally efficient algorithms, e.g., proposed in our paper. However, due to the integrality of $x$ and $y$, there may be a gap between lower bounds and upper bounds of the primal problem. To further reduce the duality gap, Fisher et al. (1997) introduce a branch-and-bound method and use the variable splitting approach to control the lower bounds; but using a new branch-and-bound method decreases the computational efficiency of our algorithm dramatically.

In our research, in order to address the similar concerns, we also apply set partitioning approach to enumerate all possible passengers' service patterns. To define passengers' service patterns, we utilize the path representation for the Traveling Salesman Problem (TSP) suggested by Bellman (1962) and Held and Karp (1962). Service pattern $j$ is defined as a vector consists of $|P|$ number of elements ($P$ is the set of passengers). Note that $p$th element of pattern $j$ is representative of passenger $p$'s service status. The service status of passenger $p$ is chosen from the set $\{0,1,2\}$, where 0 means passenger $p$ is still waiting to be picked up, 1 means passenger $p$ has been picked up but the service has not been completed, and 2 means passenger $p$ has been successfully delivered. Let $c_{vj}$ denote the travel cost of pattern $j$ traversed by vehicle $v$. Moreover, assume that $\alpha_{vpj}$ is a binary constant which equals 1 if pattern $j$ traversed by vehicle $v$ includes passenger $p$, and 0 otherwise. $z_{vj}$ is also a binary variable equals 1 if pattern $j$ is used by vehicle $v$, and 0 otherwise. Thus, we will have:

$$min \sum_{v \in V} \sum_{j \in J} c_{vj} z_{vj} \qquad (21)$$
s.t.
$$\sum_{v \in V} \sum_{j \in J} \alpha_{vpj} z_{vj} = 1 \qquad \forall p \quad (22)$$
$$\sum_{j \in J} z_{vj} = 1 \qquad \forall v \quad (23)$$
$$z_{vj} \in \{0,1\} \qquad \forall v,j \quad (24)$$

In this formulation, objective function (21) minimizes the total travel cost. Constraints (22) guarantee that each passenger is served exactly once. Constraints (23) ensure that each vehicle selects only one pattern. Constraints (24) define that the decision variables are binary. In order to assess the solution optimality of our developed algorithm on instances proposed by Ropke and Cordeau (2009) and avoid the computational challenges, two scenarios have been examined. First, we test our algorithm for the single-vehicle cases to avoid the complexity of assigning different vehicles to multiple passengers. In this case, it is obvious that we do not need to apply the set partitioning method mentioned above. Second, we test our algorithm on the small subsets of their instances for the multiple-vehicle cases with a limited number of transportation requests, so that all possible combinations of passenger-to-vehicle assignment patterns $\alpha_{vpj}$ can be enumerated, and then solved in the set partitioning problem defined above. In both sets of scenarios, the exact solutions are obtained for all the restricted master problems in the test data sets.

## 7 Conclusions

A new generation of transportation network companies uses mobile-phone-based platforms to seamlessly connect drivers to passengers from different origins to different destinations with specific, preferred departure or arrival times. Many relevant practical aspects need to be carefully formulated for real-world planning/dispatching system deployment, such as time-dependent link travel times on large-scale regional transportation networks, and tight vehicle capacity and passenger service time window constraints.

By reformulating the PDPTW through space–time networks to consider time window requirements, our proposed approach can not only solve the vehicle routing and scheduling problem directly in large-scale transportation networks with time-dependent congestion, but also avoid the complex procedure to eliminate any sub-tour possibly existing in the optimal solution for many existing formulations. By further introducing virtual vehicle constructs, the proposed approach can fully incorporate the full set of interacting factors between passenger demand and limited vehicle capacity in this model to derive feasible solutions and practically important system-wide cost-benefit estimates for each request through a sub-gradient-based pricing method. This joint optimization and pricing procedure can assist transportation network service providers to quantify the operating costs of spatially and temporally distributed trip requests.

On a large-scale regional network, the capacity impact of optimized passenger-to-vehicle matching results can be further evaluated in mesoscopic dynamic traffic simulation packages such as DTALite (Zhou and Taylor, 2014). Future work will concentrate on the development of the model for the following cases: (*i*) Passengers may desire different ride-sharing capacities (i.e. a passenger may desire to share his ride with up to only one passenger, whereas the other passenger may have no restriction about the number of passengers which share their ride with him). (*ii*) A passenger may desire to be or not to be served by a particular vehicle. (*iii*) A transportation request could contain a group of passengers who have the same origin, while they may or may not have the same destination. Alternatively, a transportation request could contain a group of passengers who have the same destination, while they may or may not have the same origin. In this case, we are interested in adding dummy nodes corresponding to passengers' origins and destinations more wisely and efficiently. In addition, in our future research, a comprehensive branch-and-bound algorithm should be included in our solution framework to fully address the complexity of assigning different vehicles to multiple passengers.

**Acknowledgments**

This paper is partially supported by National Science Foundation – United States under Grant No. CMMI 1538105 "Collaborative Research: Improving Spatial Observability of Dynamic Traffic Systems through Active Mobile Sensor Networks and Crowdsourced Data" , a U.S. University Transportation Center project titled "scheduling and managing self-driving cars for enhanced transportation system mobility and safety", and project RCS2015K006 sponsored by State Key Laboratory of Rail Traffic Control and Safety, Beijing Jiaotong University, China. We thank Prof. Hani S. Mahmassani at Northwestern University, Prof. Pitu Mirchandani at Arizona State University, and Dr. Lingyun Meng at Beijing Jiaotong University for their valuable comments. We would also like to thank our colleague, Jeffrey Taylor at the University of Utah for his help throughout the project. The work presented in this paper remains the sole responsibility of the authors.

**Appendix A. Description of the PDPTW in the origin–destination network**

Cordeau (2006) formulated the PDPTW on a network that is built based on demand request nodes and the links are defined as direct connections between pickup and delivery nodes (without explicitly considering transportation links or paths). For a systematic comparison, the following notation is adapted from Cordeau (2006).

**Table A.1**
Sets, indices and parameters used in Cordeau (2006) for the PDPTW.

| Symbol | Definition |
|---|---|
| $n$ | Number of passengers |
| $P$ | Set of passengers' pickup nodes. $P = \{1, \ldots, n\}$ |
| $D$ | Set of passengers' delivery nodes. $D = \{n+1, \ldots, 2n\}$ |
| $0$ | Node representative of origin depot |
| $2n+1$ | Node representative of destination depot |
| $N$ | Set of passengers' pickup and drop-off nodes and vehicles' depots. $N = \{P, D, \{0, 2n+1\}\}$ |
| $A$ | Set of arcs |
| $G$ | Directed graph $G = (N, A)$ |
| $i$ | Passenger $i$'s pickup node |
| $n+i$ | Passenger $i$'s delivery node |
| $q_i$ | Load at node $i$, ($i \in N$) |
| $d_i$ | Service duration at node $i$, ($i \in N$) |
| $e_i$ | Earliest time at which service is allowed to start at node $i$, ($i \in N$) |
| $l_i$ | Latest time at which service is allowed to start at node $i$, ($i \in N$) |
| $(i,j)$ | Index of arc between adjacent nodes $i$ and $j$ |
| $c_{ij}$ | Routing cost of arc $(i,j)$ |
| $t_{ij}$ | Travel time of arc $(i,j)$ |
| $V$ | Set of vehicles |
| $v$ | Vehicle index |
| $Q_v$ | Capacity of vehicle $v$ |
| $T_v$ | Maximal duration of vehicle $v$'s route |
| $L$ | Maximum ride time of a passenger |

Note that $q_0 = q_{2n+1} = 0$, $q_i \geq 0$ for $(i = 1, \ldots, n)$, and $q_i = -q_{i-n}$ $(i = n+1, \ldots, 2n)$, and service duration $d_i \geq 0$ and $d_0 = d_{2n+1} = 0$. Time window $[e_i, l_i]$ is also specified either for the pickup node or for the drop-off node of a request, but not for both. The arc set is also defined as $A = \{(i,j): (i = 0, j \in P) \text{ or } (i \in P \cup D, j \in P \cup D, i \neq j, i \neq n+j) \text{ or } (i \in D, j = 2n+1)\}$. The model uses three-index variables $x_{ij}^v$ being equal to 1 if and only if vehicle $v$ travels from node $i$ to node $j$. Let $B_i^v$ be the time at which vehicle $v$ begins servicing node $i$ and $Q_i^v$ be the load of vehicle $v$ upon departing from node $i$. Finally, for each passenger $i$, let $L_i^v$ be the ride time of passenger $i$ on vehicle $v$. The PDPTW can be formulated as follows:

$$Min \sum_{v \in V} \sum_{i \in N} \sum_{j \in N} c_{ij}^v x_{ij}^v \tag{A.1}$$

s.t.

$$\sum_{v \in V} \sum_{j \in N} x_{ij}^v = 1 \qquad \forall i \in P \tag{A.2}$$

$$\sum_{j \in N} x_{ij}^v - \sum_{j \in N} x_{n+i,j}^v = 0 \qquad \forall i \in P, v \in V \tag{A.3}$$

$$\sum_{j \in N} x_{0j}^v = 1 \qquad \forall v \in V \tag{A.4}$$

$$\sum_{j \in N} x_{ji}^v - \sum_{j \in N} x_{ij}^v = 0 \qquad \forall i \in P \cup D, v \in V \tag{A.5}$$

$$\sum_{i \in N} x_{i,2n+1}^v = 1 \qquad \forall v \in V \tag{A.6}$$

$$x_{ij}^v (B_i^v + d_i + t_{ij}) \leq B_j^v \qquad \forall i \in N, j \in N, v \in V \tag{A.7}$$

$$x_{ij}^v (Q_i^v + q_j) \leq Q_j^v \qquad \forall i \in N, j \in N, v \in V \tag{A.8}$$

$$L_i^v = B_{n+i}^v - (B_i^v + d_i) \qquad \forall i \in P, v \in V \tag{A.9}$$

$$B_{2n+1}^v - B_0^v \leq T_v \qquad \forall v \in V \quad (A.10)$$
$$e_i \leq B_i^v \leq l_i \qquad \forall i \in N, v \in V \quad (A.11)$$
$$t_{i,n+i} \leq L_i^v \leq L \qquad \forall i \in P, v \in V \quad (A.12)$$
$$max\{0, q_i\} \leq Q_i^v \leq min\{Q_v, Q_v + q_i\} \qquad \forall i \in N, v \in V \quad (A.13)$$
$$x_{ij}^v \in \{0,1\} \qquad \forall i \in N, j \in N, v \in V \quad (A.14)$$

The objective function (A.1) minimizes the total routing cost. (A.2) guarantees that each passenger is definitely picked up. (A.2) and (A.3) ensure that each passenger's origin and destination are visited exactly once by the same vehicle. (A.4) expresses that each vehicle $v$ starts its route from the origin depot. (A.5) ensures the flow balance on each node. (A.6) expresses that each vehicle $v$ ends its route at the destination depot. (A.7) and (A.8) ensure the validity of the time and load variables. (A.9) defines each passenger's ride time. (A.10) to (A.13) impose maximal duration of each route, time windows, the ride time of each passenger, and capacity constraints, respectively. Since the non-negativity of the ride time of each passenger guarantees that node $i$ is visited before node $n + i$, (A.12) also functions as precedence constraints.

**List of Tables**

**List of Figures**